\documentclass[11pt,a4paper]{article}
\usepackage{amsmath, amssymb}
\usepackage{latexsym, color}
\usepackage{epsfig}

\numberwithin{equation}{section}

\newtheorem{thm}{Theorem}[section]

\newtheorem{lem}[thm]{Lemma}

\def\nm{\noalign{\medskip}}

\newcommand{\qed}{\hfill \ensuremath{\square}}
\newcommand{\ds}{\displaystyle}
\newcommand{\pf}{\noindent {\sl Proof}. \ }
\newcommand{\p}{\partial}

\newcommand{\norm}[1]{\| #1 \|}

\newcommand{\eqnref}[1]{(\ref {#1})}

\newcommand{\Rbb}{\mathbb{R}}

\newcommand{\Bcal}{\mathcal{B}}

\newcommand{\Ncal}{\mathcal{N}}

\newcommand{\Rcal}{\mathcal{R}}

\def\Ba{{\bf a}}

\def\Be{{\bf e}}

\def\Bi{{\bf i}}
\def\Bj{{\bf j}}

\def\Bo{{\bf o}}
\def\Bp{{\bf p}}

\newcommand{\Ga}{\alpha}
\newcommand{\Gb}{\beta}
\newcommand{\Gd}{\delta}
\newcommand{\Ge}{\epsilon}

\newcommand{\Gf}{\phi}
\newcommand{\Gvf}{\varphi}
\newcommand{\Gg}{\gamma}

\newcommand{\Gk}{\kappa}

\newcommand{\Gl}{\lambda}
\newcommand{\Gn}{\eta}
\newcommand{\Gm}{\mu}

\newcommand{\Gt}{\theta}

\newcommand{\Gr}{\rho}

\newcommand{\Gs}{\sigma}

\newcommand{\Gj}{\tau}

\newcommand{\Gz}{\zeta}
\newcommand{\GD}{\Delta}

\newcommand{\GG}{\Gamma}

\newcommand{\GO}{\Omega}

\newcommand{\beq}{\begin{equation}}
\newcommand{\eeq}{\end{equation}}
\newcommand{\ol}{\overline}

\begin{document}

\title{Precise estimates of the field excited by an emitter in presence of closely located inclusions of a bow-tie shape\thanks{\footnotesize This work was
supported by NRF grants No. 2015R1D1A1A01059212, 2016R1A2B4011304 and 2017R1A4A1014735.}}

\author{Hyeonbae Kang\thanks{\footnotesize Department of Mathematics and Institute of Applied Mathematics, Inha University, Incheon
22212, S. Korea (hbkang@inha.ac.kr).} \and KiHyun Yun\thanks{\footnotesize Department of Mathematics, Hankuk University of Foreign Studies, Yongin-si, Gyeonggi-do 17035, S. Korea (kihyun.yun@gmail.com).}}

\maketitle

\begin{abstract}
This paper studies in a quantitatively precise manner the field enhancement due to presence of an emitter of the dipole type near the bow-tie structure of perfectly conducting inclusions in the two-dimensional space. We put special emphasis on field enhancement near vertices of the bow-tie structure, and derive upper and lower bounds of the gradient blow-up there. All three different kinds of symmetries are considered by varying locations and directions of the emitter, and a different estimate is derived for each case.
\end{abstract}



\section{Introduction}

In presence of closely located inclusions, there may occur field concentration in the narrow region between inclusions, which may be regarded as strain in the context of elasticity and field enhancement in the electro-static context. Regarding this phenomenon, the following mathematical model has been widely studied:
\beq\label{main}
\begin{cases}
\GD u = 0 \quad\mbox{in } \Rbb^2 \setminus \overline {(\GO_{1} \cup \GO_{2})}, \\
\ds u = {c}_j  \quad\mbox{on }\p \GO_{j}, \ \ j=1,2,  \\
\ds \int_{\p \GO_{j}} \p_{\nu} u ds =0, \quad j=1,2,\\
\ds u (X) - \Ba \cdot X = O(|X|^{-1}) \quad\mbox{as } |X| \rightarrow \infty,
\end{cases}
\eeq
where $\GO_1$ and $\GO_2$ denote the inclusions, and the gradient of the solution, $\nabla u$, represents either the electrical field or the strain. Here, ${\nu}$ is the inward-pointing, i.e., pointing toward $\GO_{1} \cup \GO_{2}$, unit normal vector field on $\p (\GO_1 \cup \GO_2)$. The second condition that $u$ takes a constant value $c_j$ on $\p\GO_j$ indicates that $\GO_j$ are perfectly conducting (the conductivity being $\infty$). The constant values $c_j$ are determined by the third condition in \eqnref{main}. The function $\Ba \cdot X$ is the linear potential which is the solution without inclusions.

Let
$$
\Ge:= \mbox{dist} (\GO_1, \GO_2).
$$
Then, the problem regarding the model \eqnref{main} is to quantitatively understand the behavior of $\nabla u$ as $\Ge$ tends to zero. This problem was raised in relation to analysis of stress in composites with stiff inclusions \cite{bab, keller}. It is also related to the effective medium theory of densely packed perfect conductors \cite{FK-CPAM-73, Keller-JAP-63}. There has been significant progress on this problem in last two decades or so: The optimal blow-up rates of $|\nabla u|$ have been derived in two dimensions \cite{AKL-MA-05, Y}, and in three dimensions \cite{BLY-ARMA-09}. The singular behavior of $\nabla u$ is characterized asymptotically near the narrow region in between two inclusions in \cite{ACKLY-ARMA-13, KLY-MA-15, KLY-JMPA-13, KLY-SIAP-14}. It is worth mentioning that the estimate for the gradient blow-up was extended to cases of insulating inclusions \cite{AKL-MA-05, BLY-CPDE-10, Yun-JDE-16}. It is also related to the spectral properties of the Neumann-Poincar\'e operator corresponding to two inclusions \cite{ACKLY-ARMA-13, BT, BT2}. Quite recently, the gradient blow-up in presence of inclusions of the bow-tie shape has been investigated in \cite{KY}. It is shown in quantitatively precise manner that the corner singularity of the solution is amplified near the vertices by interaction between closely located inclusions. We also refer to References in \cite{KY} for a more comprehensive list of references on this development.

The gradient blow-up in the model \eqnref{main} occurs passively, namely, solely by presence of inclusions and the potential $\Ba \cdot X$. However, in some contexts, the gradient blow-up, or the field enhancement, is actively created. For example, to achieve field enhancement on the bow-tie antenna, a dipole type emitter can be added to the structure (see, for example, \cite{PBFLN}). The purpose of this paper is to investigate the field enhancement when an emitter is located near closely situated inclusions. Therefore, the problem that we consider in this paper can be formulated as
\beq\label{def_gov_1}
\begin{cases}
\GD u = \Ba\cdot \nabla \Gd_{\Ge {\Be}}  \quad&\mbox{in } \Rbb^2 \setminus \ol {(\GO_{1} \cup \GO_{2})}, \\
\ds u =c_j  \quad&\mbox{on }\p \GO_{j},  \\
\ds \int_{\p\GO_{j} } \p_{\nu} u ds  =0, ~&j=1,2,\\
\ds u (X) = O\left( |X|^{-1}\right)~&\mbox{as }|X|\rightarrow\infty.
\end{cases}
\eeq
Here, $\Ba\cdot \nabla \Gd_{\Ge {\Be}}$ represents the emitter of dipole type: the unit vector $\Ba$ indicates the direction of the dipole and $\Ge {\Be}$ its location. We assume that $\Be = (0,p)$ for some $p$ with $|p| \le C$ for some $C$, say $C=1$, so that the emitter located on the $x_2$-axis. As we mention below, the inclusions are symmetric with respect to the $x_2$-axis, and so the emitter is located $\Ge$-order away from the boundaries of the inclusions. This paper deals with the case when $\GO_1$ and $\GO_2$ are bow-tie shape separated by the distance $\Ge$. In a companion paper \cite{KY2}, the circular inclusion case, as a typical example of inclusions with smooth boundaries, is dealt with.

To describe inclusions of bow-tie shape precisely, let $\GG_1$ and $\GG_2$ be the open cones in the left-half and right-half spaces with the vertices at $S_1=(-1/2,0)$ and $S_2=(1/2,0)$, respectively, and let $\Ga$ be the common aperture angle of $\GG_1$ and $\GG_2$ at $S_1$ and $S_2$ so that
\beq\label{GGj}
\GG_j= \{ Y=(y_1,y_2) ~:~ |y_2| < (-1)^j \tan(\Ga/2) (y_1-(-1)^j1/2) \}, \quad j=1,2.
\eeq
Two inclusions $\GO_1$ and $\GO_2$ for the bow-tie structure can be defined locally as translates of $\GG_1$ and $\GG_2$, namely,
\beq\label{GOone}
\GO_1 \cap B_\Gm= (\GG_1+L_\Ge) \cap B_\Gm \quad\mbox{and}\quad \GO_2 \cap B_\Gm= (\GG_2+R_\Ge) \cap B_\Gm,
\eeq
where
\beq\label{LGe}
L_\Ge:= ((-\Ge +1)/2,0) \quad\mbox{and}\quad R_\Ge:= ((\Ge-1)/2,0).
\eeq
Here $B_\Gm$ denotes the open disk of radius $\Gm$ centered at the origin $\Bo=(0,0)$. We assume that $\Gm>1$ just for ease of presentation of results. We further assume that $\GO_1$ and $\GO_2$ have smooth boundaries except at vertices $V_1:=(-\Ge/2,0)$ and $V_2:=(\Ge/2,0)$, and they are symmetric with respect to both $x_1$- and $x_2$-axes. One can easily see that the following relation holds as well:
\beq\label{GOone-2}
\Ge (\GG_1 \cup \GG_2) \cap B_\Gm = (\GO_1 \cup \GO_2) \cap B_\Gm.
\eeq
It is worthwhile to emphasize that the shapes of $\GO_1$ and $\GO_2$ are independent of $\Ge$.

For any point $\Bp$ in the plane let
\beq\label{NcalBp}
\Ncal_{\Bp}(X) = \frac 1 {2\pi} \log |X- \Bp |, \quad X \neq \Bp.
\eeq
Then it holds that $\GD (\Ba \cdot \nabla \Ncal_\Bp )= \Ba \cdot \nabla \Gd_\Bp$. So one can expect that the gradient of the solution $u$ to \eqnref{def_gov_1} may have singularity of size $|X-\Ge\Be|^{-2}$. Since the distance between vertices and $\Ge\Be$ is of order $\Ge$, the singularity caused by the emitter is expected to be of size $\Ge^{-2}$ near vertices. The main purpose of this paper is to look into the question whether there occurs enhancement of $\nabla u$ beyond $\Ge^{-2}$, especially near vertices.

We will deal with three different cases by varying the direction $\Ba$ and the location $\Ge\Be$ of the emitter: (i) $\Ba=(1,0)$, (ii) $\Ba=(0,1)$ and $\Be=(0,0)$, (iii) $\Ba=(0,1)$ and $\Be \neq (0,0)$. Each case exhibits a different symmetry: (i) the problem \eqnref{def_gov_1} is skew-symmetric with respect to the $x_2$-axis, (ii) skew-symmetric with respect to the $x_1$-axis and symmetric with respect to the $x_2$-axis, (iii) symmetric with respect to the $x_2$-axis only. We will show that in cases (i) and (iii) the field is enhanced near vertices. For instance, we obtain the following estimate for the solution $u$ to \eqnref{def_gov_1} near the vertices $V_j$:
\beq\label{nearV}
\left| \nabla u(X) \right| \simeq \frac{1}{\Ge^{1+\Gb}} \frac{1}{|X-V_j|^{1-\Gb}},
\eeq
where $\Gb$ is a number determined by the aperture angle (see \eqnref{Gb}). It is helpful to mention that $|X-V_j|^{-1+\Gb}$ is the corner singularity of the elliptic problem found by Kontratiev \cite{Kondra-TMMS-67} (see also \cite{Grisvard-book, KMR-book}). We also prove that there is no field enhancement in the case (ii), namely, $|\nabla u(X)| \lesssim \Ge^{-2}$. Here and throughout this paper, the expression $A\lesssim B$ implies that there exists a positive constant $C$ (independent of $\Ge$) such that $A \leq CB$, and $A \simeq B$ implies that both $A \lesssim B$ and $B \lesssim A$ hold. We refer to Theorems \ref{main1}, \ref{main2}, \ref{main3-1} and \ref{main3-2} for precise statements of main results.

It is instructive to compare \eqnref{nearV} with the estimate of the gradient of the solution to \eqnref{main} obtained in \cite{KY}: if $\widetilde{u}$ is the solution to \eqnref{main}, then the following holds near $V_j$:
\beq\label{blow-up1}
\left| \nabla \widetilde{u} (X) \right| \simeq \frac{1}{\Ge^\Gb |\log \Ge|} \frac{1}{|X-V_j|^{1-\Gb}}.
\eeq
Both \eqnref{nearV} and \eqnref{blow-up1} show that the elliptic corner singularity $|X-V_j|^{-1+\Gb}$ is amplified. But in \eqnref{nearV} it is due to presence of emitter, while in \eqnref{blow-up1} it is due to interaction between two inclusions.

The field enhancement shown in \eqnref{nearV} is mainly due to the corner singularity, not the interaction between two inclusion as pointed out in the sentence right after Lemma \ref{estutwo}. In the case of circular inclusions, the field enhancement is due to the interaction between two inclusions, and its magnitude is increased by the factor $\Ge^{-1/2}$ (see \cite{KY2}). It is quite interesting to observe that the factor $\Ge^{-1/2}$ is the same as that for the problem \eqnref{main} of strictly convex inclusions with smooth boundaries.

This paper is organized as follows. In the next section we review some preliminary results obtained in \cite{KY} and derive some new estimates. In section \ref{sec:decom}, we presents a decomposition of the solution to \eqnref{def_gov_1} which constitutes the basic framework of investigation in the following sections. Three sections to follow are to prove main results related to three cases mentioned above.

\section{Auxiliary functions and their estimates}

In this section we review some properties of auxiliary functions introduced in \cite{KY}. We also derive some new estimates for such functions, which will be used in later sections.

Let $\GG_j$ be the open cone defined before in \eqnref{GGj}. We choose and fix a point $P_j$ on $\p\GG_j \cap \{ (y_1,y_2): y_2>0 \}$ for $j=1,2$, and define $\Gt_j(Y)$ and $r_j(Y)$ for $Y \in \Rbb^2 \setminus \ol {\GG_j}$ by
\beq\label{angle}
\Gt_j(Y):= \mbox{the angle between $\overrightarrow {S_j P_j}$ and $\overrightarrow{S_j Y}$},
\eeq
and
\beq\label{rldef}
r_j(Y) :=  |Y-S_j|.
\eeq
The coordinates $(r_j,\Gt_j)$ may be regarded as the polar coordinates with respect to $S_j$.

Let
\beq\label{Gb}
\Gb := \frac{\pi}{2\pi - \Ga},
\eeq
where $\Ga$ is the common aperture angle of $\GG_j$. Define the functions $\Bcal_j(Y)$ for $j=1,2$ by
\beq
\Bcal_j(Y) = r_j (Y)^{\Gb} \sin (\Gb \Gt_j (Y)), \quad Y \in \Rbb^2 \setminus \ol{\GG_j} .
\eeq
The function $\Bcal_j$ is the singular part of the solution to the elliptic problem when the domain has a corner of the angle $2\pi-\Ga$ (see, for example, \cite{Grisvard-book}). The gradient $\nabla \Bcal_j$ has a singularity of order $\Gb-1$ at the vertex $S_j$. In fact, we have
\beq\label{nablaBcal}
|\nabla \Bcal_j(Y)| = \Gb |Y-S_j|^{\Gb-1}, \quad j=1,2.
\eeq
The following estimate is proved in \cite[Proposition 2.2]{KY}:
\beq\label{Bone-Btwo}
\Gb \sum_{j=1}^2 |Y-S_j|^{\Gb-1} \simeq |\nabla (\Bcal_1-\Bcal_2)(Y)|, \quad Y \in \Rbb^2 \setminus (\GG_1 \cup \GG_2).
\eeq

Now let
\beq\label{Pi}
\Pi:= \Rbb^2 \setminus \ol {\GG_1 \cup \GG_2} \quad\mbox{and}\quad \Pi^+:= \{ (y_1,y_2)\in \Pi, ~y_2>0 \}.
\eeq
Let $Q$ be the intersection point of two straight lines containing the line segments $\ol{P_1 S_1}$ and $\ol{P_2 S_2}$, where $P_j$ is the point chosen before. One can see that $Q$ is given by
$Q= \left(0, - \tan ({\Ga}/2))/2 \right)$.
We then define, for $Y \in \overline {\Pi^+}$,
\beq\label{phi_1st_deft}
\Gf(Y) := \mbox{the acute angle between } \overrightarrow {Q P_1}\mbox{ and }\overrightarrow {QY}.
\eeq
It is obvious that
$$
0 < \Gf < \pi - \Ga = \frac{\pi}{\Gg} \quad\mbox{in } \Pi^+,
$$
where
\beq\label{Ggdef}
\Gg := \frac {\pi} {\pi - \Ga} > 1.
\eeq
We also define $\Gr=\Gr(Y)$ to be the distance between $Y$ and $Q$, namely, $ \Gr(Y):= |Y-Q|.$ Thus, $(\Gr, \Gf)$ is the polar coordinate system with respect to $Q$ in $\Pi^+$. One can see from the definition of $\Gf$ that
\beq\label{nablaGf}
|\nabla \Gf(Y)| = \frac{1}{|Y-Q|}, \quad Y  \in \overline{\Pi^+}.
\eeq

Let $w$ be the solution to
\beq\label{wcondi}
\begin{cases}
\GD w = 0 \quad\mbox{in }  \Pi^+, \\
w= 0  \quad\mbox{on } \p\Pi^+ \setminus \left[-\frac 1 2 , \frac 1 2 \right]\times \{0\} ,  \\
\nm
\p_2 w =  - \p_2 \Gf \quad\mbox{on } \left[-\frac 1 2 , \frac 1 2 \right]\times \{0\} ,\\
\nm
\ds \int_{\Pi^+} |\nabla w|^2 dA < \infty ,
\end{cases}
\eeq
where $\p_j$ denotes the partial derivative with respect to the $x_j$-variable for $j=1,2$. Note that $\Gf$ and $w$ are defined in $\Pi^+$. We extend them to $\Pi$ as symmetric functions with respect to the $y_1$-axis, namely,
\beq
\Gf(y_1,y_2) = \Gf(y_1,-y_2), \quad w(y_1,y_2) = w (y_1,-y_2). \label{ext_Gf_w}
\eeq

It is proved in \cite[Lemma 3.1 (ii)]{KY} that
$$
|\nabla w (Y) | \lesssim  \frac{1}{|Y-Q|^{\Gg + 1 }}, \quad Y \in \Pi^+ \setminus  { B_2} ,
$$
from which one can immediately see that the following holds:
\beq\label{nablaw}
|\nabla w (Y) | \lesssim  \frac{1}{|Y|^{\Gg + 1 }}, \quad Y \in \Pi \setminus  { B_2}.
\eeq
Here and throughout this paper $B_r(P)$ denotes the open disk of radius $r$ centered at $P$, and if $P=\Bo$, the origin, then we simply denote it by $B_r$. It is also proved in \cite[Lemma 3.2]{KY} that
there is a constant $a$, $0<a< 2^{\Gb+2}/\Gg$, such that
\beq\label{phiwBcal}
|\nabla (\Gf + w)(Y) - a \nabla (\Bcal_1 - \Bcal_2 ) (Y)| \leq C , \quad Y \in \Pi.
\eeq

Let $q$ be the solution to
\beq\label{qeqn}
\begin{cases}
\GD q = 0 \quad\mbox{in }  \Rbb^2 \setminus \ol {(\GO_{1} \cup \GO_{2})}, \\
\ds q = \Gl_j  \quad\mbox{on } \p \GO_{j}, \ \ j=1,2, \\
\ds \int_{\p \GO_{1}} \p_{\nu} q ds = -\ds \int_{\p \GO_{2}} \p_{\nu} q ds =-1, \\
\ds q (X)  = O(|X|^{-1}) \quad\mbox{as } |X| \rightarrow \infty  .
\end{cases}
\eeq
Here, $\Gl_1$ and $\Gl_2$ are constants determined by the third and fourth conditions in \eqnref{qeqn}, and they depend on $\Ge$. Note that
\beq
q |_{\p \GO_1}=\Gl_1 < 0 < \Gl_2 = q |_{\p \GO_2}
\eeq
and
\beq\label{hopf2}
(-1)^j \p_{\nu}q  > 0 \quad\mbox{on } \p \GO_j \setminus \{V_j\} , \quad j=1,2,
\eeq
which is a consequence of Hopf's lemma. It is worth mentioning that we are using the notation $X = (x_1,x_2)$ to denote points in $\Rbb^2 \setminus (\GO_1 \cup \GO_2)$, and $Y = (y_1, y_2)$ for points in the scaled region, namely, $\Rbb^2 \setminus (\GG_1 \cup \GG_2)$.

The function $q$ was introduced and played a crucial role in \cite{KY}, and so does in this paper. It is proved in \cite[Lemma 4.3]{KY} that there exists $\Ge_0 >0$ such that
\beq\label{lemma43}
q|_{\p \GO_2} - q|_{\p \GO_1} \simeq \frac 1 {|\log \Ge|}
\eeq
for all $\Ge \le \Ge_0$.

We now define the function $l$ by
\beq\label{decomp_q_phi_w}
q (X)=  (q|_{\p \GO_2} - q|_{\p \GO_1}) \left( \frac {\Gg} {\pi}(\Gf + w ) (\Ge^{-1} X) + l (X) \right), \quad X \in B_\Gm \setminus \ol{ \GO_1 \cup \GO_2} .
\eeq
Then $l$ is a harmonic function and satisfies
\beq
l|_{(\p \GO_1) \cap B_1} =  l|_{(\p \GO_2) \cap B_1} = - \frac 1 2 ,
\eeq
because $\lambda_1 = -\lambda_2$, $(\Gf + w ) (\Ge^{-1} X)=0$ on $(\p \GO_1) \cap B_1$ and $(\Gf + w ) (\Ge^{-1} X)=\pi/\Gg$ on $(\p \GO_2) \cap B_1$.

The following lemma collects estimates for $q$ and $l$ which are used in later sections.

\begin{lem}\label{lem_q_behavior}
There exists $\Ge_0>0$ such that the following holds for all $\Ge \le \Ge_0$:
\begin{itemize}
\item[(i)] There is a constant $\Gd$ independent of $\Ge$ such that
\beq\label{nablaq1}
| \nabla q (X)| \lesssim q|_{\p \GO_2} - q|_{\p \GO_1},  \quad X \in \Rbb^2 \setminus \ol {B_{\Gd} \cup \GO_1 \cup \GO_2},
\eeq and\beq\label{nablaq2}
|\nabla q (X)| \lesssim \frac 1 {\Ge^\Gb |\log \Ge|} \sum_{j=1}^2 \frac{1}{|X-V_j|^{1-\Gb}}, \quad X \in    \p ( \Omega_1 \cup  \Omega_2) \cap B_{\delta }  \setminus \{V_1,V_2\} .
\eeq
\item[(ii)] There is a constant $\Gj$ independent of $\Ge$ such that
\beq\label{pnuq}
(-1)^j \p_{\nu} q (X) \gtrsim  \frac 1 {\Ge |\log \Ge|}, \quad X \in \p \GO_j \cap B_{\Gj \Ge} (V_j) \setminus \{V_j\}.
\eeq
\item[(iii)] It holds that
\beq\label{nablal}
|\nabla l (X)| \lesssim \sum_{j=1}^2 \frac{1}{|X-V_j|^{1-\Gb}}
\eeq
for $X \in B_{1} \setminus ( {\GO_1 \cup \GO_2}  \cup \{V_1,V_2\})$.
\end{itemize}
\end{lem}

\pf
By the maximum principle, we have $q|_{\p \GO_1} \leq q \leq q|_{\p \GO_2}$, and
$$|q| \leq  q|_{\p \GO_2} - q|_{\p \GO_1}.$$

Fix $\Gd < 1 $. There are two simply connected domains, say $U_1$ and $U_2$, with smooth boundaries such that
$$
\GO_1 \setminus B_{\Gd/2} \subset U_1 \subset \GO_1 \quad \mbox{and} \quad \GO_2 \setminus B_{\Gd/2} \subset U_2 \subset \GO_2.
$$
Here we assume that $\Ge$ is so small that the vertices $V_1$ and $V_2$ belong to $B_{\Gd / 4}$. By the Riemann mapping theorem, there are two conformal mappings $\Psi_{i}: \Rbb^2 \setminus {\overline {B_1}} \rightarrow \Rbb^2 \setminus {\overline {U_j}}$ ($i=1,2$). Since $\p U_1$ and $\p U_2$ are smooth, there is a constant $C >0$ such that
\beq\label{CC}
C \leq |\nabla \Psi_{i} (Z)| \leq C^{-1}
\eeq
for all $Z \in \Rbb^2\setminus \overline{B_1}$ and for $i=1,2$. (See, for example, equation (29) of \cite{KLY-MA-15} for a proof of this fact.) Here the constant $C$ may be taken to be independent of $\Ge$. In fact, let $U_i^0$ be the domain $U_i$ when $\Ge=\Ge_0$, and $\Psi_i^0$ be the corresponding conformal mapping. Then we may take $U_i$ corresponding to $\Ge$ to be the translate of $U_i^0$, namely, $U_i= U_i^0 + \frac{(-1)^i}{2}(\Ge-\Ge_0)$, and corresponding conformal mapping can be taken to be $\Psi_i =\Psi_i^0 + \frac{(-1)^i}{2}(\Ge-\Ge_0)$. So, the constant $C$ in \eqnref{CC} can be taken to be independent of $\Ge$. In particular, one can show using \eqnref{CC} that there is a constant $C_*$ independent of $\Ge$ such that
\beq\label{C*}
C_* |Z_1-Z_2| \le |\Psi_1(Z_1)-\Psi_1(Z_2)| \le C_*^{-1} |Z_1-Z_2|
\eeq
for all $Z_1$ and $Z_2$ in $\Rbb^2 \setminus {\overline {B_1}}$.

Let $V$ be the left half of $\Rbb^2 \setminus \ol {B_{\Gd} \cup \GO_1 \cup \GO_2}$, i.e.,
$$
V:= ((-\infty,0] \times \Rbb) \setminus \overline{B_\Gd \cup \GO_1}.
$$
We claim that there is $\Gd_0$ such that for all $Z \in \Psi_{1}^{-1} (V)$, $q \circ\Psi_{1}$ can be extended into $B_{\Gd_0} (Z)$ as a harmonic function and the modulus of the extended function is bounded by $2(q|_{\p \GO_2} - q|_{\p \GO_1})$.

To prove the claim, we suppose that $Z \in \Psi_{1}^{-1} (V)$ and $\mbox{dist}(Z, \p \Psi_{1}^{-1} (V)) < C_* \Gd/2$ where $C_*$ is the constant appearing in \eqnref{C*}. It is helpful to mention that for the other case, namely the case when $\mbox{dist}(Z, \p \Psi_{1}^{-1} (V)) \ge C_* \Gd/2$, the claim holds trivially by taking $\Gd_0= C_* \Gd/2$. Let $X:= \Psi_1(Z)$. Then, by \eqnref{C*}, we see that $\mbox{dist}(X, \p V) < \Gd/2$. So there are three cases to happen : (i) $\mbox{dist}(X, \p\GO_1 \setminus B_{\Gd}) < \Gd/2$, (ii) $\mbox{dist}(X, \p\GO_1 \setminus B_{\Gd}) \ge \Gd/2$ and $\mbox{dist}(X, B_{\Gd}) < \Gd/2$, (iii) $X$ is within $\Gd/2$-distance from the $y$-axis. In the third case the claim holds trivially. In the second case, we have $B_{\Gd/2}(X) \subset \Rbb^2 \setminus \overline{ B_{\Gd/2} \cup \GO_1 \cup \GO_2}$. So, $q \circ\Psi_{1}$ is harmonic in $\Psi_1^{-1}(B_{\Gd/2}(X))$, and by \eqnref{C*} $\Psi_1^{-1}(B_{\Gd/2}(X))$ contains $B_{C_*^{-1} \Gd/2}(Z)$.
To prove the claim in the first case, we observe that $q \circ\Psi_{1}$ is constant on $\Psi_{1}^{-1} (\p \GO_1 \setminus B_{\Gd/2})$, which is a subset of $\p B_1$. So there is $\Gd_1$ such that for each $Y \in \Psi_{1}^{-1} (\p \GO_1 \setminus B_{\Gd})$ $q \circ\Psi_{1}$ is extended to $B_{\Gd_1}(Y)$ by reflection. Because the extended function is defined by reflection, the modulus of the extended function is bounded by $2(q|_{\p \GO_2} - q|_{\p \GO_1})$. Thus the claim holds by taking $\Gd_0$ to be the minimum of $C_*^{-1} \Gd/2$ and $\Gd_1$.

Here we invoke a standard gradient estimate for harmonic functions: if $h$ is a harmonic function defined in an open set containing a closed ball $\overline {B_r (X_0)}$, then
\beq\label{standard}
|\nabla h(X_0)| \leq  \frac 2 {r}  \norm {h}_{L^{\infty} (B_r (X_0))}.
\eeq
This inequality will be repeatedly used in this paper.

If $Z \in \Psi_{1}^{-1} (V)$, then $q \circ\Psi_{1}$ can be extended into $B_{\Gd_0} (Z)$ as a harmonic function and the modulus of the extended function is bounded by $2(q|_{\p \GO_2} - q|_{\p \GO_1})$. It then follows from \eqnref{standard} that
$$
|\nabla (q (\Psi_{1\Gd}) )(Z) | \leq \frac 8 {\Gd_0} (q|_{\p \GO_2} - q|_{\p \GO_1}) .
$$
We then use \eqnref{CC} to infer that
$$
|\nabla q  (X)) | \leq \frac {8\sqrt 2} {C \Gd_0} (q|_{\p \GO_2} - q|_{\p \GO_1})
$$
for all $X \in V$.

So far we proved that \eqnref{nablaq1} holds on the left half of $\Rbb^2 \setminus \ol {B_{\Gd} \cup \GO_1 \cup \GO_2}$.
By the symmetry of $\GO_1 \cup \GO_2$, the same inequality also holds in the right half of $\Rbb^2 \setminus \ol {B_{\Gd} \cup \GO_1 \cup \GO_2}$.

Let us now prove \eqnref{nablal}. It is proved in \cite[Lemma 4.1]{KY} that \eqnref{nablal} holds for $X \in B_{\Gd_1} \setminus  ( {\GO_1 \cup \GO_2} \cup \{V_1,V_2\})$ for some $\Gd_1$. By taking smaller $\Gd$ in \eqnref{nablaq1} or $\Gd_1$ if necessary, we may assume that $\Gd_1=\Gd$. We now prove \eqnref{nablal} for $X \in B_1 \setminus (B_{\Gd} \cup   { \GO_1 \cup \GO_2})$ for sufficiently small $\Ge$ so that $V_1$ and $V_2$ do not belong to $B_1 \setminus (B_{\Gd} \cup   { \GO_1 \cup \GO_2})$. If $X \in \Rbb^2 \setminus ({B_{\Gd} \cup \GO_1 \cup \GO_2})$, it follows from \eqnref{nablaGf} that
\beq\label{nablaphi}
|\nabla (\Gf (\Ge^{-1} X))| \lesssim \frac 1 {\Ge}  \frac 1 { |(x_1,|x_2|)/\Ge-Q|}  \lesssim 1.
\eeq
It then follows from the definition \eqnref{decomp_q_phi_w} of $l$, together with \eqnref{nablaw}, \eqnref{lemma43}, \eqref{nablaq1} and \eqnref{nablaphi}, that
\begin{align*}
|\nabla l (X)|& \lesssim \left( q|_{\p \GO_2} - q|_{\p \GO_1} \right)^{-1} | \nabla q (X)| + \frac{\Gg}{\pi} \left( |\nabla (\Gf (\Ge^{-1} X))| + |\nabla (w (\Ge^{-1} X))| \right)  \\
& \lesssim 1 \\
& \lesssim \frac 1 {|X-V_1|^{1-\beta}} + \frac 1 {|X-V_2|^{1-\beta}}
\end{align*}
for any $X \in B_{1} \setminus \left(B_{\Gd} \cup \ol {\GO_1 \cup \GO_2}\right)$. So we have \eqnref{nablal}.

One can see from \eqnref{nablaBcal}, \eqnref{phiwBcal}, \eqnref{decomp_q_phi_w} and \eqnref{nablal} that the following inequality holds:
\begin{align*}
|\nabla q (X)| & \le  (q|_{\p \GO_2} - q|_{\p \GO_1}) \left( \frac{\Gg}{\pi} \big| \Ge^{-1} \nabla (\Gf + w ) (\Ge^{-1} X) \big| + |\nabla l (X)| \right) \\
& \lesssim (q|_{\p \GO_2} - q|_{\p \GO_1}) \left( \sum_{j=1}^2 \Ge^{-1} \big| \nabla \Bcal_j (\Ge^{-1} X) \big| + |\nabla l (X)| + \Ge^{-1}\right)  \\
& \lesssim (q|_{\p \GO_2} - q|_{\p \GO_1}) \Ge^{-\beta} \sum_{j=1}^2 \frac{1}{|X-V_j|^{1-\Gb}}
\end{align*} for any $X $ with $|X| < 2\Ge$ on  $\p ( {\GO_1 \cup \GO_2} ) \setminus \{V_1,V_2\} $. On the other hand, it follows from \eqref{nablaGf}, \eqref{nablaw}, \eqnref{decomp_q_phi_w} and \eqnref{nablal} that
\begin{align*}
|\nabla q (X)| & \le  (q|_{\p \GO_2} - q|_{\p \GO_1}) \left( \frac{\Gg}{\pi} \big| \Ge^{-1} \nabla (\Gf + w ) (\Ge^{-1} X) \big| + |\nabla l (X)| \right) \\
& \lesssim (q|_{\p \GO_2} - q|_{\p \GO_1}) \left( \Ge^{-1} \frac {\Ge}{ |X|} +  \Ge^{-1}  \left(\frac {\Ge}{ |X|} \right)^{1+\gamma} + \frac 1 {|X|^{1-\beta}}  \right)  \\
& \lesssim (q|_{\p \GO_2} - q|_{\p \GO_1})   \Ge^{-1}  \left(\frac {\Ge}{ |X|} \right)^{1-\beta}   \\
& \lesssim (q|_{\p \GO_2} - q|_{\p \GO_1}) \Ge^{-\beta} \sum_{j=1}^2 \frac{1}{|X-V_j|^{1-\Gb}}
\end{align*} for any $X $ with $1 > |X| \geq  2\Ge$ on  $\p ( {\GO_1 \cup \GO_2} )  \setminus  \{V_1,V_2\} $. So, \eqnref{nablaq2} follows from  \eqnref{lemma43}.

We now prove \eqnref{pnuq}, only for $j=2$. The case when $j=1$ can be treated in the exactly same manner. According to \eqnref{phiwBcal}, we have
$$
\p_{\nu} \left( (\Gf + w)( \Ge^{-1} X) \right) \geq  a \Ge^{-1} \p_{\nu} (\Bcal_1 - \Bcal_2 )( \Ge^{-1} X)  - C \Ge^{-1}.
$$
From \eqnref{nablaBcal}, one can find a proper $\tau_0$ such that
$$
\p_\nu \Bcal_2(Y) = \Gb |Y-S_2|^{\Gb-1},
$$
and
$$
|\nu \cdot \nabla \Bcal_1(Y)| \le \Gb |Y-S_1|^{\Gb-1} \le \frac{1}{2} \p_\nu \Bcal_2(Y)
$$
for all $Y \in \p \GG_2 \cap B_{\tau_0} (S_2) \setminus \{S_2\}$. So, we have for all $X \in \p \GO_2 \cap B_{\tau_0 \Ge} (V_2) \setminus \{V_2\}$
$$
\p_{\nu} \left( (\Gf + w)( \Ge^{-1} X) \right) \geq  \Ge^{-1} \left( \frac {a\Gb} 2  \frac {\Ge^{1-\Gb}} {| X-V_2|^{1-\Gb}}  - C \right).
$$
We then infer that there is a constant $\Gj < \Gj_0$ such that
\beq\label{pnu}
\p_{\nu} \left( (\Gf + w)( \Ge^{-1} X) \right) \gtrsim  \Ge^{-\beta} \frac {1}{|X-V_2|^{1-\beta}}
\eeq
for all $X \in \p \GO_2 \cap B_{\Gj\Ge} (V_2) \setminus \{V_2\}$. Now \eqnref{pnuq} follows from \eqnref{decomp_q_phi_w}. Indeed, \eqnref{decomp_q_phi_w} yields
$$
\p_\nu q (X) =  (q|_{\p \GO_2} - q|_{\p \GO_1}) \left( \frac{\Gg}{\pi} \p_{\nu} \left( (\Gf + w)( \Ge^{-1} X) \right) + \p_\nu l (X) \right).
$$
Hence, we see from  \eqnref{lemma43}, \eqnref{nablal} and \eqnref{pnu} that there is $\Ge_0>0$ such that  \eqnref{pnuq} holds for $\Ge< \Ge
_0$. This completes the proof. \qed

\section{Decomposition of the solution}\label{sec:decom}

Let $u$ be the solution to \eqnref{def_gov_1}, and define the function $\Gs$ by
\beq\label{def_c_a_e_e_2}
u = \frac{u|_{\p \GO_2}- u|_{\p \GO_1}}{q|_{\p \GO_2}- q|_{\p \GO_1}} q + \Gs  \quad \mbox{in } \Rbb^2 \setminus ({\GO_1 \cup \GO_2 \cup \{\Ge \Be \}}).
\eeq
Then in view of \eqnref{decomp_q_phi_w}, we may write
$$
u(X) = u_1(X) + \Rcal_1(X) + \Gs(X), \quad X \in B_{\mu} \setminus \ol{ \GO_1 \cup \GO_2},
$$
where $\mu$ is the number introduced in \eqref{GOone-2} which we assume to be larger than 1,
\beq\label{uone}
u_1(X) := \frac{\Gg (u|_{\p \GO_2}- u|_{\p \GO_1} ) }{\pi} (\Gf + w ) (\Ge^{-1} X),
\eeq
and
\beq\label{Rone}
\Rcal_1(X):=(u|_{\p \GO_2}- u|_{\p \GO_1} ) l(X).
\eeq

Observe that $\GD \Gs = \Ba \cdot \nabla  \Gd_{\Ge {\Be}}$ in $\Rbb^2 \setminus \ol {(\GO_{1} \cup \GO_{2})}$ and $\Gs$ attains the same constant values on $\p\GO_1$ and $\p \GO_2$, namely,
\beq\label{Gsbdry}
\Gs |_{\p \GO_1} = \Gs|_{\p \GO_2}= \frac{u|_{\p \GO_2}+ u|_{\p \GO_1}}{2} .
\eeq
Let $v$ be the solution to
\beq\label{def_v}
\begin{cases}
\GD v = 0 \quad \mbox{in } \Rbb^2 \setminus \ol {(\GO_{1} \cup \GO_{2})}, \\
\ds v  = {\Ba} \cdot \nabla  \Ncal_{\Ge \Be }  \quad\mbox{on }\p \GO_{1} \cup \p \GO_{2},  \\
\ds \int_{\Rbb^2 \setminus \ol{ \left(\GO_{1} \cup \GO_{2} \right ) }} |\nabla v|^2  dA < \infty,
\end{cases}
\eeq
where $\Ncal_\Bp$ is defined by \eqnref{NcalBp}.
Then we have
\beq\label{sigma_N_v}
\nabla \Gs = \nabla \left( {\Ba} \cdot \nabla  \Ncal_{\Ge \Be } - v \right) \quad\mbox{in }\Rbb^2 \setminus \ol {(\GO_{1} \cup \GO_{2})}.
\eeq
Having in mind the relation
\beq\label{scaling}
\nabla  \Ncal_{\Ge \Be}(X)  = \Ge^{-1} (\nabla  \Ncal_{\Be})\left({\Ge}^{-1}X\right),
\eeq
we define $V$ to be the solution to \beq\label{def_V}
\begin{cases}
\GD V = 0 \quad \mbox{in } \Rbb^2 \setminus \ol {(\Gamma_{1} \cup \Gamma_{2})}, \\
\ds V = {\Ba} \cdot \nabla \Ncal_{ \Be }  \quad \mbox{on }\p \Gamma_{1} \cup \p \Gamma_{2},  \\
\ds \int_{\Rbb^2 \setminus \ol{ \left(\Gamma_{1} \cup \Gamma_{2} \right ) }} |\nabla V|^2 dA < \infty.
\end{cases}
\eeq

Let
\beq\label{defGvf}
\Gvf(Y) := {\Ba} \cdot \nabla  \Ncal_{\Be }(Y) - V(Y), \quad Y \in \Rbb^2 \setminus (\Gamma_{1} \cup \Gamma_{2}),
\eeq
and
\beq\label{utwo}
u_2(X) := \Ge^{-1} \Gvf (\Ge^{-1} X), \quad X \in B_{\mu} \setminus (\GO_1 \cup \GO_2).
\eeq
Then, we have
\beq\label{Gsutwo}
\nabla \Gs(X)= \nabla u_2(X) + \nabla \Rcal_2(X),
\eeq
where
\beq\label{Rtwo}
\Rcal_2(X):= -v(X)+ \Ge^{-1} V(\Ge^{-1} X).
\eeq

In summary, we have
\beq\label{udecom}
\nabla u= \nabla u_1 + \nabla u_2 + \nabla\Rcal_1 + \nabla\Rcal_2,
\eeq
where $u_j$ and $\Rcal_j$ ($j=1,2$) are defined by \eqnref{uone}, \eqnref{Rone}, \eqnref{utwo}, and \eqnref{Rtwo}, respectively.
In sections to follow, we show with precise estimates that $\nabla u_1$ and $\nabla u_2$ constitute the major terms characterizing the behavior of $\nabla u$, while $\nabla\Rcal_1$ and $\nabla\Rcal_2$ are error terms. Here we emphasize that $u_1$ is zero if the potential difference is zero as one can see from \eqnref{uone}, while $u_2$ does not have potential difference as \eqnref{Gsbdry} and \eqnref{Gsutwo} show.

\section{Case 1: ${\Ba} = {\Bi }$}

In this section we obtain optimal estimates of $\nabla u$ where $u$ is the solution to \eqnref{def_gov_1} when the direction $\Ba$ of the dipole is given by
\beq\label{BaBi}
\Ba=\Bi=(1,0).
\eeq
We first emphasize that in this case the equation in \eqref{def_gov_1} takes the form $\GD u = \p_1 \Gd_{\Ge {\Be}}$ and $\p_1 \Gd_{\Ge {\Be}}$ is odd in the $x_1$-variable. By considering symmetry of the problem, one immediately see that $u|_{\p\GO_2}=- u|_{\p\GO_1}$. So, the blow-up of $\nabla u$ in this case is caused not only by the presence of the emitter but also by the potential gap $u|_{\p\GO_2} - u|_{\p\GO_1}$ (and presence of the corners).

The following is the main result of this section.

\begin{thm}\label{main1}
Let $u$ be the solution to \eqref{def_gov_1} under the assumption \eqnref{BaBi}, and define the region $R$ by
\beq\label{regionR}
R := B_1 \setminus \overline { \GO_1 \cup \GO_2 \cup \{\Ge {\bf e}\}}.
\eeq
There exist $c_0 \in (0,1/2)$ and  $\Ge_0>0$ such that the following estimates hold for all $\Ge \le \Ge_0$:
\begin{itemize}
\item[(i)] For all $X \in  ( B_{c_0 \Ge}(V_1)  \cup B_{c_0 \Ge }(V_2) ) \cap R$,
\beq\label{unear}
| \nabla u (X)  | \simeq \frac 1 {\Ge^{1+\Gb}} \sum_{j=1}^2 \frac{1}{|X - V_j|^{1-\Gb} },
\eeq
where $\Gb$ is the number defined by \eqnref{Gb}.
\item[(ii)] For all $X \in R$ with $c_0^{-1} \Ge |\log \Ge| < |X| < c_0$,
\beq\label{uaway}
| \nabla u  (X) | \simeq \frac 1 {\Ge | \log \Ge|} \frac{1}{|X|}.
\eeq

\item[(iii)]  For all $X \in R \setminus (  B_{c_0 \Ge}(V_1) \cup B_{c_0 \Ge }(V_2))$,
\end{itemize}
\beq\label{uaway-upper-bd}
| \nabla u  (X) | \lesssim  \frac{1}{|X - \Ge {\bf e} |^2} + \frac 1 {\Ge | \log \Ge|} \frac{1}{|X| + \Ge}.
\eeq
\end{thm}

The estimate \eqnref{unear} clearly shows that the field $\nabla u$ is enhanced near the vertex $V_j$, from $\Ge^{-2}$ (due to presence of dipole type emitter) to $\Ge^{-1-\Gb} |X - V_j|^{-1+\Gb}$.

To prove this theorem we obtain the following lemmas, which are about estimates of four functions appearing in the decomposition \eqnref{udecom}. Proofs of these lemmas are given in subsections to follow.

\begin{lem}\label{estuone}
Let $R$ be the region defined by \eqnref{regionR}. There is a constant $\Ge_0 >0 $  such that the following estimates for $u_1$ defined by \eqnref{uone} hold for all $\Ge \le \Ge_0$:
\begin{itemize}
\item[(i)] There exists $\Gk_1 \in (0,1/2)$ such that
\beq\label{uonenear}
| \nabla u_1 (X) | \simeq \frac{1}{\Ge^{1+\Gb} |\log \Ge|} \sum_{j=1}^2 \frac{1}{|X - V_j|^{1-\Gb}}
\eeq for all $X \in ( B_{\Gk_1 \Ge}(V_1)  \cup B_{\Gk_1 \Ge }(V_2) ) \cap R$.
\item[(ii)]  There exists  $\Gk_2 >2 $ such that
\beq\label{uoneaway}
| \nabla u_1 (X) | \simeq \frac{1}{\Ge |\log \Ge|} \frac{1}{|X|}
\eeq for all $X \in R$ with $|X|> \Gk_2 \Ge $.
\item[(iii)]  For any $\Gk_3 \in (0,1/2) $, it holds that
\beq\label{uoneaway_upper}
| \nabla u_1 (X) | \lesssim \frac{1}{\Ge |\log \Ge|} \frac{1}{|X| + \Ge}
\eeq for all $X \in R \setminus (  B_{\Gk_3 \Ge}(V_1) \cup B_{\Gk_3 \Ge }(V_2)) $.
\end{itemize}
\end{lem}

\begin{lem}\label{estutwo}
Let $R$ be the region defined by \eqnref{regionR}. There is  $\Ge_0$ such that the following estimates for $u_2$ defined by \eqnref{utwo} hold for all $\Ge \le \Ge_0$:
\begin{itemize}
\item[(i)] There exists $\Gk_1 \in (0,1/2)$ such that
\beq\label{utwonear}
\left| \nabla u_2(X) \right| \simeq \frac{1}{\Ge^{1+\Gb}} \sum_{j=1}^2 \frac{1}{|X-V_j|^{1-\Gb}}
\eeq
for all $X \in ( B_{\Gk_1 \Ge}(V_1)  \cup B_{\Gk_1 \Ge }(V_2) ) \cap R$.

\item[(ii)]  For any $\Gk_2 \in (0,1/2) $, it holds that
\beq\label{utwoaway}
|\nabla u_2(X)| \lesssim |X - \Ge {\bf e}|^{-2}
\eeq for all $X \in R \setminus (  B_{\Gk_2 \Ge}(V_1) \cup B_{\Gk_2 \Ge }(V_2))  $.
\end{itemize}
\end{lem}

We emphasize that blow-up in \eqnref{uonenear} is weaker than that in \eqnref{utwonear}. Since $u_1$ carries information on the potential difference of the solution $u$ to \eqnref{def_gov_1} while $u_2$ does not, \eqnref{uonenear} and \eqnref{utwonear} suggest that the field enhancement due to presence of corner is stronger than that due to potential difference. In this regard we add in subsection \ref{subsec:single} a proof showing that blow-up estimate \eqnref{utwonear} is valid even when there is a single inclusion with a corner, not a bow-tie structure.

\begin{lem}\label{estRone}
There exists $\Ge_0>0$ such that
\beq\label{nablaRone}
| \nabla \Rcal_1 (X)  | \lesssim \frac 1 {\Ge |\log \Ge|} \sum_{j=1}^2 \frac 1 {|X - V_j|^{1-\Gb}}
\eeq
for all $X \in B_{1} \setminus \ol {\GO_1 \cup \GO_2}$ and $\Ge \le \Ge_0$.
\end{lem}

\begin{lem}\label{estRtwo}
There exists  $\Ge_0>0$ such that
\beq\label{nablaRtwo}
\left| \nabla \Rcal_2(X) \right| \lesssim \sum_{j=1}^2 \frac{1}{|X-V_j|^{1-\Gb}}
\eeq
for all $X \in B_1 \setminus \ol {\GO_{1}\cup \GO_{2}}$ and $\Ge \le \Ge_0$.
\end{lem}

Theorem \ref{main1} is a consequence of these lemmas. In fact, \eqnref{unear} is an easy consequence of \eqnref{uonenear}, \eqnref{utwonear}, \eqnref{nablaRone} and \eqnref{nablaRtwo}. We then see from \eqnref{utwoaway} that
$$
|\nabla u_2(X)| \le \frac{c_0 C}{\Ge |\log \Ge|} \frac{1}{|X|}
$$
for some constant $C$ provided that $c_0^{-1} \Ge |\log \Ge| < |X|$. So, \eqnref{uaway} follows from \eqnref{uoneaway}, \eqref{utwoaway}, \eqnref{nablaRone} and \eqnref{nablaRtwo} if $c_0$ is small enough. Finally \eqref{uaway-upper-bd} follows from \eqref{uoneaway_upper},  \eqref{utwoaway}, \eqref{nablaRone} and \eqref{nablaRtwo}.

\subsection{Proofs of Lemmas \ref{estuone} and  \ref{estRone} }
We begin with the following lemma.

\begin{lem}\label{lem:potential_difference}
Let $u$ be the solution to \eqref{def_gov_1} under the assumption \eqnref{BaBi}. There exists $\Ge_0>0$ such that the following holds for all $\Ge \le \Ge_0$:
\beq
u|_{\p \GO_2}- u|_{\p \GO_1} \simeq \frac 1 {\Ge \left| \log \Ge \right|}.
\eeq
\end{lem}

\pf
We first note that
$$
\int_{\p \GO_1 \cup \p \GO_2} \p_{\nu} ( u - \p_1 \Ncal_{\Ge {\Be}}) ds = \int_{\Rbb^2 \setminus (\GO_1 \cup \GO_2)} \GD ( u - \p_1 \Ncal_{\Ge {\Be}}) dA=0.
$$
So, it follows from the definition \eqref{qeqn} of $q$ and the divergence theorem that
\begin{align}
u|_{\p \GO_2}- u|_{\p \GO_1} & = \int_{\p \GO_1 \cup \p \GO_2} u \p_{\nu} q ds \notag \\
& = \int_{\p \GO_1 \cup \p \GO_2} \left(\p_1 \Ncal_{\Ge {\Be}}+ ( u - \p_1 \Ncal_{\Ge {\Be}})\right) \p_{\nu} q ds \notag  \\
& = \int_{\p \GO_1 \cup \p \GO_2} (\p_1 \Ncal_{\Ge {\Be}}) \, \p_{\nu} q ds. \label{u_N_q}
\end{align}

Note that
\beq\label{pxNcal}
\p_1 \Ncal_{\Ge \Be}(X) =  \frac{1}{2\pi} \frac{x_1}{x_1^2+(x_2-\Ge p)^2},
\eeq
from which one immediately see that $\p_1 \Ncal_{\Ge \Be}(X) >0$ if $x_1>0$ and $\p_1 \Ncal_{\Ge \Be}(X) <0$ if $x_1<0$. This together with \eqnref{hopf2} and \eqnref{u_N_q} yields
$$
u|_{\p \GO_2}- u|_{\p \GO_1} \ge \int_{\p \GO_2 \cap B_{\Gj \Ge} (V_2)} \p_1 \Ncal_{\Ge {\Be}} \, \p_{\nu} q ds,
$$
where $\Gj$ is the number appeared in \eqnref{pnuq}. One can also see from \eqnref{pxNcal} that $\p_1 \Ncal_{\Ge {\Be}} (X) \simeq \Ge^{-1}$ for all $X \in \p \GO_2 \cap B_{\Gj \Ge} (V_2)$. It then follows from \eqnref{pnuq} that
$$
\int_{\p \GO_2 \cap B_{\Gj \Ge} (V_2)} \p_1 \Ncal_{\Ge {\Be}} \, \p_{\nu} q ds \gtrsim \frac 1 {\Ge |\log \Ge|},
$$
and hence
$$
u|_{\p \GO_2}- u|_{\p \GO_1} \gtrsim \frac 1 {\Ge |\log \Ge|}.
$$

To prove the opposite inequality, let $\Gd$ be the number appeared in Lemma \ref{lem_q_behavior} (i), and write
$$
u|_{\p \GO_2}- u|_{\p \GO_1}  = \int_{(\p \GO_1 \cup \p \GO_2) \cap B_{\Gd}} + \int_{(\p \GO_1 \cup \p \GO_2) \setminus B_{\Gd}} (\p_1 \Ncal_{\Ge {\Be}}) \, \p_{\nu} q ds  =:I + II.
$$
Since $|\p_1 \Ncal_{\Ge {\Be}} (X)| \lesssim |X|^{-1}$
for any $X \in \p \GO_1 \cup \p \GO_2$, we obtain using \eqnref{nablaq2}
\begin{align*}
|I| &\le \frac 1 {\Ge^\Gb |\log \Ge|} \int_{(\p \GO_1 \cup \p \GO_2) \cap B_{\Gd}}
\sum_{j=1}^2 \frac{1}{|X-V_j|^{1-\Gb} |X|} ds\\
&\lesssim \frac 1 {\Ge^\Gb |\log \Ge|} \sum_{j=1}^2 \int_{(\p \GO_j) \cap B_{\Gd}}
\frac{1}{|X-V_j|^{1-\Gb} (\Ge + |X-V_j|)} ds \lesssim \frac 1 {\Ge |\log \Ge|}.
\end{align*}
On the other hand, using \eqnref{lemma43} and \eqnref{nablaq1} we obtain
$$
|II| \lesssim \frac 1 { |\log \Ge|} \int_{(\p \GO_1 \cup \p \GO_2) \setminus B_{\Gd}} \frac 1 {|X|} ds \lesssim  \frac 1 { |\log \Ge|}.
$$
Thus we have
$$
u|_{\p \GO_2}- u|_{\p \GO_1}  \lesssim \frac 1 {\Ge |\log \Ge|},
$$
and the proof is complete. \qed

In view of the definition \eqnref{Rone} of $\Rcal_1$, Lemma \ref{estRone} immediately follows from \eqnref{nablal} and Lemma \ref{lem:potential_difference}.

\medskip

\noindent{\sl Proof of Lemma \ref{estuone}}.
We first see from \eqnref{nablaGf} and \eqnref{nablaw} that there is a constant $C_1$ such that
$$
\frac{1}{|Y|} - \frac{C_1}{|Y|^{\Gg + 1 }} \lesssim |\nabla (\Gf+ w) (Y)| \lesssim \frac{1}{|Y|} + \frac{1}{|Y|^{\Gg + 1 }}
$$
for all $Y \in \Pi$ with $|Y| >2$. So, if $|X| > 2\Ge$, then
$$
\frac{\Ge}{|X|} - C_1 \left( \frac{\Ge}{|X|} \right)^{\Gg + 1 } \lesssim |\nabla (\Gf+ w) (\Ge^{-1} X)| \lesssim \frac{\Ge}{|X|} + \left( \frac{\Ge}{|X|} \right)^{\Gg + 1 }.
$$
So, there is $\Gk_2 >2$ such that if $|X| \ge \Gk_2 \Ge$, then
\beq\label{uoneest1}
|\nabla (\Gf+ w) (\Ge^{-1} X)| \simeq \frac{\Ge}{|X|}.
\eeq

On the other hand, we infer from \eqnref{Bone-Btwo} and \eqnref{phiwBcal} that there is $\Gk_1$ such that
\beq\label{uoneest2}
|\nabla (\Gf + w)(\Ge^{-1}X)| \simeq |\nabla (\Bcal_1 - \Bcal_2 )(\Ge^{-1}X)| \simeq \sum_{j=1}^2 \frac{\Ge^{1-\Gb}}{|X - V_j|^{1-\Gb}},
\eeq
for $X$ satisfying $|X-V_j|\le \Gk_1 \Ge$ ($j=1,2$).

Now one can see from \eqnref{uone} and Lemma \ref{lem:potential_difference} that \eqnref{uonenear} and \eqnref{uoneaway} for all $\Ge < \Ge_0$ follow from  \eqnref{uoneest2} and \eqnref{uoneest1}, respectively. Here, $\Ge_0$ is as given in Lemma \ref{lem:potential_difference}.

The function $u_1$, defined in $B_{\mu} \setminus \overline {\GO_1 \cup \GO_2}$ by \eqref{uone}, where $\mu$ is the constant in \eqref{GOone}, is positive and bounded by $u|_{\p \GO_2} - u|_{\p \GO_1}$. For any $\Gk_3 \in (0,1/2) $, there exists a constant $C_2$ such that for all $X \in R \setminus ( B_{\Gk_3 \Ge}(V_1) \cup B_{\Gk_3 \Ge }(V_2) )$ and $\Ge < \Ge_0$,
$$
B_{C_2 (|X|+\Ge)} (X) \subset B_{\mu} \setminus \{V_1,V_2\}.
$$
Since $u_1$ is constant on $\p \GO_i \cap B_{\mu}$ for $i = 1,2$, $u_1$ can be extended into $B_{C_2 (|X|+\Ge)} (X)$ as a harmonic function whose modulus is bounded by $2 (u|_{\p \GO_2} - u|_{\p \GO_1})$.  Thus, \eqref{standard} yields
$$
| \nabla u_1 (X) | \lesssim (u|_{\p \GO_2} - u|_{\p \GO_1}) \frac{1}{|X| + \Ge} .
$$
Hence, we obtain \eqref {uoneaway_upper}  by Lemma \ref{lem:potential_difference}, and the proof is complete.
\qed

\subsection{Proof of Lemma \ref{estutwo}}

We first prove the following lemma.

\begin{lem}\label{lem_bounds}
Let $v$ and $V$ be the solutions of \eqref{def_v} and \eqref{def_V}, respectively, when ${\Ba} = {\Bi }$.
The following holds.
\begin{itemize}
\item [(i)] $\p_1 \Ncal_{ \Ge \Be } < v < 0$  in $\left(\Rbb^- \times  \Rbb\right) \setminus \ol {\GO_{1}}$,
 \item [(ii)] $\p_1 \Ncal_{ \Ge \Be } > v > 0$ in $\left(\Rbb^+ \times  \Rbb\right) \setminus \ol {\GO_{2}}$,
\item [(iii)] $\p_1 \Ncal_{  \Be } < V <  0$ in $\left(\Rbb^- \times  \Rbb\right) \setminus \ol {\Gamma_{1}}$,
\item [(iv)] $\p_1 \Ncal_{  \Be } > V > 0$ in $\left(\Rbb^+ \times  \Rbb\right) \setminus \ol {\Gamma_{2}}$.
\end{itemize}
\end{lem}

\pf
We only prove (iv); the others can be proved similarly.

We first note that
\beq
\p_1 \Ncal_{\Be }(Y) = \frac{y_1}{2\pi |Y-{\Be}|^2},
\eeq
in particular, $\p_1 \Ncal_{\Be}$ is odd with respect to the $y_1$-axis. By the symmetry of the inclusions, we infer that $V$ is also odd with respect to the $y_1$-axis. In particular, $V(Y)=0$ if $y_1=0$. Since $V=\p_1 \Ncal_{\Be} >0$ on $\p\GG_{2}$, $V>0$ on $R_0:=(\Rbb^+ \times \Rbb) \setminus \ol {\GG_{2}}$ by the maximum principle.

Since $V(Y)=0$ if $y_1=0$, $V$ is of the form $y_1F(Y)$ for some $F$. Thus there are $M>0$ and $r_0>0$ such that
$$
M y_1> V(Y) > 0
$$
for $|Y-{\Be} | < r_0$ with $y_1 >0$. If $r<r_0$ is sufficiently small so that $2\pi r^2 < M^{-1}$, then
$$
\p_1 \Ncal_{\Be}(Y) > M y_1 > V(Y) >0
$$
for all $Y$ with $|Y-\Be|=r$. So, the function $\p_1  \Ncal_{\Be} - V$ is harmonic in $R_0\setminus B_r(\Be)$ and positive on $\p(R_0\setminus B_r(\Be))$. Since
$$
\int_{R_0\setminus B_r(\Be)} |\nabla (\p_1  \Ncal_{\Be} - V)|^2 dA < \infty ,
$$
we infer from the maximum principle that $\p_1  \Ncal_{\Be} - V >0$ in $R_0\setminus B_r(\Be)$. Since $r$ is arbitrary, we arrive at (iv).
\qed

\begin{lem}\label{lem:B}
Let $\Gvf$ be the function defined by \eqnref{defGvf}. There exists a constant $c \in (0,  1 / 2)$ such that  \beq\label{eqn:B}
\left| \nabla \Gvf(Y) \right| \simeq \sum_{j=1}^2 \frac{1}{|Y-S_j|^{1-\Gb}}
\eeq
for all $Y \in ( B_{c }(S_1) \setminus \ol {\GG_{1}} ) \cup ( B_{c }(S_2) \setminus \ol {\GG_{2}} )$.
\end{lem}

\pf
Recall that $\Gvf= \p_1 \Ncal_{\Be} - V$. Let $(r_2, \Gt_2)$ be the polar coordinates with respect to $S_2$ as defined in \eqref{angle} and \eqref {rldef}.  Since $\Gvf$ vanishes on $\p\GG_2$, it admits the Fourier series expansion of the form
\beq\label{Gvffourier}
\Gvf (Y) = \sum_{n=1}^{\infty} a_n r_2 ^{n\Gb} \sin\left(n \Gb \Gt_2 \right), \quad Y \in B_{1/2}(S_2) \setminus \ol {\GG_{2}},
\eeq
where the first coefficient $a_1$ is given by
$$
a_1 = \frac{2^{2\Gb+3}}{2\pi- \Ga} \int_{\p B_{1/ 4}(S_2) \setminus \GG_{2}} \Gvf(Y) \sin \left( \Gb \Gt_2 (Y) \right) d s .
$$
By Lemma \ref{lem_bounds} (iv), $\Gvf (Y) > 0$ if $Y \in \p B_{1/4}(S_2) \setminus \ol {\GG_{2}}$.  Hence, we have $a_1>0$. We emphasize that the first order term in the expansion \eqnref{Gvffourier} is $\Bcal_2$, in other words,
\beq\label{Gvffourier2}
\Gvf (Y)- a_1 \Bcal_2(Y) = \sum_{n=2}^{\infty} a_n r_2 ^{n\Gb} \sin\left(n \Gb \Gt_2 \right).
\eeq

Let $0< c < 1/2$. We infer from \eqnref{Gvffourier} that
$$
\norm {\Gvf}_{L^{2} (B_{(1+2c)/ 4}(S_2) \setminus \GG_{2})}^2 = \frac{2\pi-\Ga}{2} \sum_{n=1}^{\infty} \frac{|a_n|^2}{2n\Gb +2} \left(\frac {1+2c} 4 \right)^{2n\Gb + 2} .
$$
On the other hand, Lemma \ref{lem_bounds} (iv) implies $0<\Gvf \le  \p_1 \Ncal_{\Be}$. We then have
$$
\norm {\Gvf}_{L^{2} (B_{(1+2c)/ 4}(S_2) \setminus \GG_{2})}^2 \le \norm {\p_1 \Ncal_{\Be}}_{L^{2} (B_{(1+2c)/ 4}(S_2) \setminus \GG_{2})}^2 \le C
$$
for some constant $C$, where the last inequality holds to be true because $\Be$ is away from $B_{(1+2c)/ 4}(S_2)$. It then follows that
\beq\label{1000}
\sum_{n=1}^{\infty} \frac{|a_n|^2}{2n\Gb +2} \left(\frac {1+2c} 4 \right)^{2n\Gb + 2} \le C.
\eeq

Suppose $Y \in B_{c}(S_2) \setminus \ol{\GG_{2}}$. Since $\Gb> \frac 1 2 $, it follows from \eqnref{Gvffourier2} that
$$
\left| \nabla \Gvf(Y) - a_1 \nabla \Bcal_2 (Y) \right| = \left| \nabla \left( \sum_{n=2}^{\infty} a_n r_2^{n\Gb} \sin(\Gb \Gt_2) \right) \right|
\leq  \sum_{n=2}^{\infty} 2n \Gb |a_n| c^{n\Gb - 1 }.
$$
It then follows from the Cauchy-Schwartz inequality and \eqnref{1000} that
\begin{align*}
&\left| \nabla \Gvf(Y) - a_1 \nabla \Bcal_2 (Y) \right| \\
& \leq \left( \sum_{n=1}^{\infty} \frac{|a_n|^2}{2n\Gb +2} \left(\frac {1+2c} 4 \right)^{2n\Gb + 2 } \right)^{1/2} \left( c^{-4} \sum_{n=2}^{\infty}  {(2n\Gb +2)}{|2n\Gb|^2} \left(\frac{4c}{1+2c}\right)^{2n\Gb +2} \right)^{1/2}
\\& \leq C_1
\end{align*}
for some constant $C_1$. Now \eqnref{eqn:B} for $Y \in B_{c }(S_2) \setminus \ol {\GG_{2}}$ follows from \eqnref{nablaBcal} if $c$ is sufficiently small.

By switching the roles of $\Bcal_j$ ($j=1,2$) in above arguments, we obtain \eqnref{eqn:B} for $X \in B_{c }(S_1) \setminus \ol {\GG_{1}}$, and the proof is complete.
\qed

\begin{lem}  For any $ k\in (0,1/2)$, it holds that
\beq\label{nablaGvf}
|\nabla \Gvf(Y)| \lesssim |Y - \Be|^{-2}
\eeq for all $Y \in \Pi \setminus \left(B_k (S_1) \cup B_k (S_2) \right)$. Here $\Pi$ is the set defined in \eqnref{Pi}.
\end{lem}
\pf
Let
$$
\Gn(Y):= \min\{ |Y-S_1|, |Y-S_2|, |Y-\Be| \}.
$$ Since $\Gvf =0$ on $\p\Pi$, for each $Y \in \Pi$ ($Y \neq S_1, S_2, \Be$), $\Gvf$ can be extended, by a reflection with respect to either $\p\GG_1$ or $\p\GG_2$,  to $B_{c\Gn(Y)}(Y)$ as a harmonic function. Here, $c$ is a constant less than 1 which depends on the aperture angle $\Ga$, but is independent of $Y$. We denote the extended function by the same notation $\Gvf$. Then, \eqnref{standard} yields
$$
|\nabla \Gvf (Y)|  \leq  \frac{2}{c\Gn(Y)} \norm {\Gvf }_{L^{\infty} (B_{c\Gn(Y)} (Y))} =\frac{2}{c\Gn(Y)} \norm {\Gvf }_{L^{\infty} (B_{c\Gn(Y)} (Y) \cap \Pi)}.
$$ Lemma \ref{lem_bounds} shows that if $Z\in B_{c\Gn(Y)}(Y) \cap \Pi$, then
$$
|\Gvf(Z)| \leq 2 |\p_1 \Ncal_{\Be }(Z)| \leq \frac{1}{\pi |Z-\Be|} \leq \frac{1}{\pi (1-c)}\frac{1}{|Y-\Be|} .$$ Thus, we have \beq
|\nabla \Gvf (Y)|  \lesssim   \frac{1}{\Gn(Y) |Y-\Be|}  .\label{new_lem_4.9_first_eq}
\eeq

Suppose that $k\in (0,1/2)$ and we show that
\beq\label{new_lem_4.9_second_eq}
\Gn(Y) \gtrsim |Y-\Be|
\eeq
for all $Y \in \Pi  \setminus \left(B_k (S_1) \cup B_k (S_2) \right)$. If $|Y| \ge 1$, then $\Gn(Y) \simeq |Y-\Be|$, since $|{\Be}| < 1$. Otherwise, $Y$ belongs to $B_1 \setminus (B_k (S_1) \cap B_k (S_2))$ so that
$$ \eta (Y) \geq \min\{ k , |Y-\Be| \}   \geq k |Y-\Be| /2 ,
$$
since $|Y-\Be| \leq 2$ and $k < 1/2$. Hence, we have \eqref{new_lem_4.9_second_eq} for all $Y \in \Pi  \setminus \left(B_k (S_1) \cup B_k (S_2) \right)$. Therefore, \eqnref{nablaGvf} follows from \eqref{new_lem_4.9_first_eq} and \eqref{new_lem_4.9_second_eq}.
\qed

Lemma \ref{estutwo} is an immediate consequence of above two lemmas.

\subsection{Proof of Lemma \ref{estRtwo}}

We first prove the following lemma.

\begin{lem}\label{lem:v-V}  Let  $\mu$ be the constant given in \eqref{GOone-2}. There exists a positive constant $\Ge_0$ such that
\beq\label{v-V}
| \Rcal_2 (X)| \leq \frac{2}{\pi}
\eeq
for all $X \in  B_{\mu}\setminus \ol {\GO_{1}\cup \GO_{2}}$ and $\Ge \le \Ge_0$.
\end{lem}

\pf
Because of \eqnref{scaling} and \eqref{Rtwo}, we have
$$
\Rcal_2(X)= -v(X) + \Ge^{-1} V \left( {\Ge}^{-1} X\right) = 0, \quad X \in (\p \GO_{1}\cup  \p \GO_{2}) \cap B_{\mu}.
$$
If $X = (x_1,x_2) \in (\p B_{\mu}) \setminus (\GO_{1}\cup \GO_{2})$, then
\eqref{scaling} and Lemma \ref{lem_bounds} yield
$$
| \Rcal_2(X) | \leq 2|\p_1 \Ncal_{\Ge {\Be}}(X) | = \frac{|x_1|}{\pi |X-\Ge {\Be}|^2} \le \frac{1}{\pi |X-\Ge {\Be}|}.
$$
Since $|X|=\mu >1$, we have $|X-\Ge {\Be}| \ge 1-\Ge p$, and hence
$$
| \Rcal_2(X) | \leq \frac 1 {\pi (1-\Ge p)}\leq \frac 2 {\pi},
$$
provided that $\Ge$ is sufficiently small to satisfy $\Ge p < 1/2$. By taking $\Ge_0=1/(2p+1)$, the maximum principle yields \eqnref{v-V}.
\qed

\medskip

\noindent{\sl Proof of Lemma \ref{estRtwo}}.
We show that
\beq\label{lem:A-main1}
\left| \nabla \Rcal_2(X) \right| \lesssim \left|\nabla \Bcal_1\left(X-L_{\Ge}\right)\right| + \left |\nabla \Bcal_2\left(X-R_{\Ge}\right)\right|
\eeq
for all $X \in B_{1} \setminus \overline {\GO_{1}\cup \GO_{2}} $ and for any $\Ge \le \Ge_0$, where $L_{\Ge}$ and $R_{\Ge}$ are given by \eqnref{LGe} and $\Ge_0$ is a small positive number. Since $S_1+L_\Ge=V_1$ and $S_2+R_\Ge=V_2$, \eqnref{nablaRtwo} follows from \eqnref{nablaBcal} and \eqnref{lem:A-main1}.

To prove \eqnref{lem:A-main1} we closely follow the argument in the proof of Proposition 2.2 of \cite{KY}. We show that
there are constants $\Ge_0$ such that the following inequality holds for all $\Ge < \Ge_0$ and for all $X \in B_{1} \setminus (\GO_1 \cup \GO_2)$:
\beq\label{bdryest}
\left| \nabla \Rcal_2(X) \right| \leq C_1 \left| \nabla \Bcal_1 (X-L_\Ge) - \nabla \Bcal_2 (X - R_\Ge) + C_2 (x_1 + 1/2 , x_2 ) \right |
\eeq
for some constants $C_1$ and $C_2$ independent of $\Ge<\Ge_0$.

Let $\psi_{1+}$ and $\psi_{1-}$ be the solutions to
\beq
\begin{cases}
\ds \GD \psi_{1+}=\GD \psi_{1-} =  0  \quad &\mbox{in } B_{\mu} \setminus \overline{\GO_1}, \\
\ds \psi_{1+}=\psi_{1-} = 0   \quad &\mbox{on } \p\GO_1 \cap B_{\mu}, \\
\ds \psi_{1+}=- \psi_{1-} = \frac 2 {\pi}   \quad &\mbox{on }  \p  B_{\mu} \setminus \overline{\GO_1}.
\end{cases} \notag
\eeq
We infer from Lemma \ref{lem:v-V} and the fact that $\Rcal_2=0$ on $\p\GO_1 \cap B_\mu$ that $\psi_{1+} - \Rcal_2$ attains its minimum on $\p\GO_1 \cap B_{\mu}$. Likewise, $\psi_{1-} - \Rcal_2$ attains its maximum on $\p\GO_1 \cap B_\mu$. So, by Hopf's lemma, we have
\beq\label{hopf1}
\p_{\nu} \psi_{1+}   < \p_{\nu} \Rcal_2 < \p_{\nu} \psi_{1-} \quad\mbox{on } \p\GO_1 \cap B_{\mu}.
\eeq

Note that there is a small $\Ge_0>0$ satisfying $\ol{B_{1}(V_1)} \subset B_{\mu}$ if $\Ge < \Ge_0$. So, $\psi_{1+}$ admits the Fourier series expansion of the form
$$
\psi_{1+} (X) = a_0 \Bcal_1 (X-L_{\Ge}) +   \sum_{n=2}^\infty a_n r_1(X-L_\Ge)^{n\Gb} \sin\big( n\Gb \Gt_1(X-L_\Ge)),
$$
where $(r_1, \Gt_1)$ is the polar coordinates with respect to $S_1$ as defined in \eqnref{angle} and \eqnref{rldef}, and the series converges for $X$ satisfying $|X-V_1| < 1+s$ for some $s>0$. Thus we have
$$
|\nabla ( \psi_{1+} (X) - a_0 \Bcal_1 (X-L_{\Ge}))|\lesssim 1,
$$
and hence
\beq\label{psi1+}
|\nabla \psi_{1+}(X)| \leq | a_0  \nabla \Bcal_1 (X-L_{\Ge}) |+  | \nabla \psi_{1+}(X) - a_0  \nabla \Bcal_1 (X-L_{\Ge}) |  \lesssim |\nabla \Bcal_1 (X-L_{\Ge})|
\eeq
for  all $X \in \overline{B_{1} (V_1)}\setminus \GO_1$. Similarly, one can prove that
\beq\label{psi1-}
|\nabla \psi_{1-}(X)| \lesssim |\nabla \Bcal_1 (X-R_{\Ge})| \quad\mbox{for all } X \in \overline{B_1 (V_1)}\setminus \GO_1.
\eeq

Since $\Bcal_1 (X-L_{\Ge}) = 0$ and $\Rcal_2=0$ on $ \p \GO_1 \cap B_1(V_1)$, it follows from \eqnref{hopf1}, \eqnref{psi1+} and \eqnref{psi1-} that for all $X \in  \p \GO_1 \cap B_1(V_1)$
\begin{align*}
\left| \nabla \Rcal_2(X) \right| &= \left| \p_\nu \Rcal_2(X) \right|
\lesssim \left|\nabla \Bcal_1 (X-L_{\Ge})\right| = \left|\p_{\nu} \Bcal_1 (X-L_{\Ge})\right| .
\end{align*}
It is proved in \cite[Proposition 2.2]{KY} that for all $X=(x_1,x_2) \in \p \GO_1 \cap B_1(V_1)$
$$
\p_{\nu} \Bcal_1 (X-L_{\Ge}) <0, \quad  -\p_{\nu} \Bcal_1 (X-R_{\Ge}) \leq 0, \quad \nu\cdot (x_1+1/2,x_2) <0.
$$
Thus we have
\begin{align*}
\left| \nabla \Rcal_2(X) \right|  & \lesssim  \left| \p_{\nu} \Bcal_1 (X-L_\Ge)\right| \\ & \leq \left| \p_{\nu} \Bcal_1 (X-L_\Ge) -  \p_{\nu} \Bcal_2 (X - R_\Ge) + C \nu \cdot \left(x_1 + 1/2, x_2\right) \right |  \\
& \leq \left| \nabla \Bcal_1 (X-L_\Ge) - \nabla \Bcal_2 (X - R_\Ge) + C \left(x_1 + 1/2, x_2\right) \right |
\end{align*}
for all $X =(x_1,x_2) \in \p {\GO_1} \cap B_1 (V_1)$ and for any $C>0$. Similarly one can show that the same inequality holds for $X \in \p {\GO_2} \cap B_1 (V_2)$.

Since $\p {\GO_j} \cap B_1 \subset \p {\GO_j} \cap B_1 (V_j)$ for $j=1,2$, we have
\beq\label{lem:A:first:eq1}
\left| \nabla \Rcal_2(X) \right| \lesssim \left| \nabla \Bcal_1 (X-L_\Ge) - \nabla \Bcal_2 (X - R_\Ge) + C \left(x_1 + 1/2, x_2\right) \right |
\eeq
for all $X \in (\p {\GO_1} \cap \p {\GO_2}) \cap B_1 $ and for any $C>0$.

On the other hand, Lemma \ref{lem:v-V} implies that $\Rcal_2$ is bounded in $( B_{1+t} \setminus B_{1-t}) \setminus (\GO_1 \cup \GO_2 )$ for some $t>0$ regardless of $\Ge< \Ge_0$. Moreover, $\Rcal_2$ is constant on $\p \GO_1$ and $\p \GO_2$. Then we can apply \eqnref{standard} and the same argument as in the proof of Lemma \ref{lem_q_behavior} to show that $|\nabla \Rcal_2|$ is bounded on $(\p B_1) \setminus (\GO_1 \cup \GO_2)$.  Since $| \nabla \Bcal_1 (X-L_\Ge) - \nabla \Bcal_2 (X - R_\Ge) |$ is bounded on $(\p B_1) \setminus (\GO_1 \cup \GO_2)$, there is a large constant $C_2$ such that
\beq\label{lem:A:first:eq2}
\left| \nabla \Rcal_2(X) \right| \lesssim \left| \nabla \Bcal_1 (X-L_\Ge) - \nabla \Bcal_2 (X - R_\Ge) + C_2 \left(x_1 + 1/2, x_2\right) \right |
\eeq
for any $X \in (\p B_1) \setminus (\GO_1 \cup \GO_2)$. This together with \eqref{lem:A:first:eq1} proves that \eqnref{lem:A:first:eq2} holds on $\p(B_1 \setminus (\GO_1 \cup \GO_2))$. We then invoke the maximum modulus theorem to show \eqnref{lem:A:first:eq2} (or \eqref{bdryest}) in $B_1 \setminus (\GO_1 \cup \GO_2)$.

We then infer from \eqnref{bdryest} that
$$
|\nabla \Rcal_2(X) | \lesssim |\nabla \Bcal_1 ( X - L_\Ge)| + |\nabla \Bcal_2 ( X - R_\Ge )| + \left| \left( x_1+1/2,x_2 \right) \right|.
$$
Since $|\nabla \Bcal_1 ( X - L_\Ge)| + |\nabla \Bcal_2 ( X - R_\Ge )| \gtrsim 1$, but $|(x_1+1/2,x_2)| \lesssim 1$ for all $X \in B_1 \setminus (\GO_1 \cup \GO_2)$, \eqref{lem:A-main1} follows. This completes the proof of \eqnref{lem:A-main1}. \qed

\subsection{A single inclusion case}\label{subsec:single}

As mentioned briefly after Lemma \ref{estutwo}, the blow-up of $|\nabla u_2|$ in \eqnref{utwonear} has nothing to do with the interaction between two inclusions. To make it clear we consider the following problem where there is a single inclusion $\GO_1$:
\beq\label{def_1_single_inclusion}
\begin{cases}
\GD u = {\Bi } \cdot \nabla  \Gd_{{\Bo}}  \quad&\mbox{in } \Rbb^2 \setminus \ol {\GO_{1}}, \\
\ds u =c_1\mbox{ (constant) } \quad&\mbox{on }\p \GO_{1},  \\
\ds u (X) = O\left( |X|^{-1}\right)~&\mbox{as }|X|\rightarrow\infty.
\end{cases}
\eeq
Here again the constant value $c_1$ prescribed on $\p\GO_1$ indicates that $\GO_1$ is a perfect conductor. However, unlike two-inclusion model \eqnref{def_gov_1}, we do not need to impose the condition
$\int_{\p \GO_1} \p_\nu  u ds = 0$ in this case to determine the constant value $c_1$. In fact, $c_1$ is determined by the third condition in \eqnref{def_1_single_inclusion} (see the beginning of the proof below).

\begin{thm}\label{thm:single}
Let $u$ be the solution to \eqnref{def_1_single_inclusion}. Then, there exists  a constant $\Gk_0 \in (0,1/2) $ such that  the following estimate holds for all $X \in B_{\Gk_0 \Ge} (V_1) \setminus \overline{\GO_1}$:
\beq\label{one-inclu_main}
| \nabla  u (X) | \simeq \frac 1 {\Ge^{1+\Gb}} \frac 1 {|X-V_1|^{1-\Gb}} .
\eeq
\end{thm}

\pf
Let $v$ be the solution to
\beq\label{def_2}
\begin{cases}
\GD v = 0 \quad&\mbox{in } \Rbb^2 \setminus \ol {\GO_{1}}, \\
\ds v =  \p_1  \Ncal_{\Bo}\quad&\mbox{on }\p \GO_{1},  \\
\ds \int_{\Rbb^2 \setminus  \GO_1} |\nabla  v|^2 dA < \infty . ~&~
\end{cases}
\eeq
Let $c_*:= \lim_{|X| \to \infty} v(X)$. The third condition in \eqnref{def_1_single_inclusion} shows that
\beq
u(X)= \p_1  \Ncal_{\Bo}(X) - v(X) + c_*, \quad X \in \Rbb^2 \setminus \ol {\GO_{1}}.
\eeq
In particular, we have $c_1=c_*$.

As before, we introduce an auxiliary function: let $V$ be the solution to
\beq\label{def_3}
\begin{cases}
\GD V = 0 \quad&\mbox{in } \Rbb^2 \setminus \ol {\GG_{1}}, \\
\ds V =   \p_1  \Ncal_{\Bo}\quad&\mbox{on }\p \GG_{1},  \\
\ds \int_{\Rbb^2 \setminus  \GG_1} |\nabla  V|^2 dA < \infty, ~&~
\end{cases}
\eeq
and define $u_2$ and $\Rcal_2$, respectively, by
\beq
u_2 (X) = \p_1  \Ncal_{\Bo}(X) - \Ge^{-1} V (\Ge^{-1} X ),
\eeq
and
\beq
\Rcal_2 (X) =  -v(X) +\Ge^{-1} V (\Ge^{-1} X ) , \quad X \in B_{\mu} \setminus \ol {\GO_{1}}.
\eeq
Then, we have
\beq
\nabla u (X)= \nabla u_2(X) + \nabla \Rcal_2 (X), \quad X \in B_{\mu} \setminus \ol {\GO_{1}} .
\eeq

First, we prove that
\beq\label{1st-ineq-one-inclusion}
|\nabla \Rcal_2(X)| \lesssim \frac 1  {\Ge} \frac 1 {|X-V_1|^{1-\Gb}} , \quad X \in B_1 \setminus \ol {\GO_{1}} .
\eeq
To do so, we have from the maximum principle that for $X \in \Rbb^2 \setminus \GO_1$
$$
|v(X)| \leq \max_{X \in \p\GO_1} \left|\p_1  \Ncal_{\Bo} (X) \right| = \left|\p_1  \Ncal_{\Bo} (V_1) \right| \lesssim \Ge^{-1}.
$$
We also have
$$
|V(Y)| \leq \max_{Y \in \p\GG_1} \left|\p_1  \Ncal_{\Bo} (Y) \right| = \left|\p_1  \Ncal_{\Bo} (S_1) \right| \lesssim 1
$$
for $Y \in \Rbb^2 \setminus \GG_1$. Thus we have
\beq\label{R2Ge-1}
|\Rcal_2 (X) | \leq | -v(X) |+|\Ge^{-1} V (\Ge^{-1} X )  | \lesssim \Ge^{-1}
\eeq
for all $X \in B_{\mu} \setminus \ol {\GO_{1}}.$

Since $\Rcal_2=0$ on $B_\Gm \cap \p\GO_1$, it admits a Fourier series expansion of the form
$$
\Rcal_2 (X) = b_1 \Bcal_1 (X-L_{\Ge}) +  \sum_{n=2} ^{\infty} b_n r_1 (X-L_\Ge)^{n\Gb} \sin (n\Gb \Gt_1 (X-L_{\Ge})), \quad X \in B_{\mu} \setminus \ol {\GO_{1}},
$$
where
$$
b_1 = \frac {2^{\Gb+2}} {2\pi -\Ga}\int_{\p B_{1/2}(V_1) \setminus \overline{\GO_1}} \Rcal_2 (X) \sin(\Gb \Gt_1 (X-L_\Ge)) ds.
$$
The estimate \eqnref{R2Ge-1} yields that $|b_1| \lesssim \Ge^{-1}$. We also have in the same way as the proof of Lemma \ref{lem:B} that
$$
\left| \nabla (\Rcal_2 (X) -b_1 \Bcal_1 (X-L_{\Ge}) ) \right| \lesssim \Ge^{-1}
$$
for all $X \in B_1\setminus \overline{\GO_1}$.  Now \eqref{1st-ineq-one-inclusion} follows. In fact, we have
$$
\left| \nabla \Rcal_2 (X) \right| \lesssim |b_1| \left| \nabla \Bcal_1 (X-L_{\Ge}) \right| + \Ge^{-1} \lesssim \frac 1  {\Ge} \frac 1 {|X-V_1|^{1-\Gb}}.
$$

Second, we prove the existence of  a constant $\Gk_0 \in (0,1/2) $ such that
\beq\label{one-inclu}
| \nabla  u_2 (X) | \simeq \frac 1 {\Ge^{1+\Gb}} \frac 1 {|X-V_1|^{1-\Gb}}
\eeq
for all $X \in B_{\Gk_0 \Ge} (V_1) \setminus \overline{\GO_1}$.  One can see that
$$
\nabla u_2(X) = \Ge^{-2 } \nabla \Gvf (\Ge^{-1} X), \quad X \in B_{\mu} \setminus \GO_1 ,
$$
where $\Gvf$ is defined by
\beq
\Gvf (Y) = \p_1  \Ncal_{\Bo}(Y) - V (Y ) , \quad Y\in \Rbb^2 \setminus \GG_1. \label{defGvf-11}
\eeq
Since $\Gvf$ vanishes on $\p \GG_1$, it admits the Fourier series expansion of the form
\beq\label{6000-1}
\Gvf (Y) = a_1 \Bcal_1 (Y) +  \sum_{n=2}^{\infty} a_{n} r_1^{n\Gb} (Y) \sin n\Gb \Gt_1 (Y)
\eeq
for $ Y \in B_{ 1/2} (S_1) \setminus \ol{\GG_1}$, where the first coefficient $a_1$ is given by
\beq\label{a_e1-1}
a_{1} = \frac{2}{(2\pi - \Ga){r_1 ^{\Gb}}} \int_{0}^{2\pi-\Ga}  \Gvf (Y(r_1,\Gt_1)) \sin \Gb \Gt_1 d\Gt_1
\eeq
for any $r_1 \in (0,1/2)$. Since $V = \p_1 \Ncal_{\Bo} < 0 $ on $\p \GG_1$ and $\p_1 \Ncal_{\Bo}|_{\p \GG_1}$ attains its minimal value $\p_1 \Ncal_{\Bo} (S_1)$ at $S_1$,
\beq\label{N_S_0_V_inequal}
\p_1 \Ncal_{\Bo} (S_1) = V(S_1) \leq  V (Y) < 0
\eeq
for all $Y \in \Rbb^2 \setminus \GG_1$. By \eqref{defGvf-11}, we have $|\Gvf (Y)| \lesssim 1$ for all $Y \in B_{1/3} (S_1) \setminus \GG_1$. It then follows from \eqref{a_e1-1} that $|a_1| \lesssim 1$. We then obtain the following estimate using the same argument as in the proof of Lemma \ref{lem:B}:
\beq
\left| \nabla (\Gvf (X) -a_1 \Bcal_1 (X) ) \right| \lesssim 1 \label{varphi-a_1B-temp}
\eeq
for all $X \in B_{1/4}(S_1)\setminus \overline{\GG_1}$.

We now show that $a_1 \neq 0$. We first have
\begin{align*}
\Gvf(Y) &=   \p_1 \Ncal_{ {\Bo}} (Y) - V (Y) \nonumber \\
&\leq   \p_1 \Ncal_{ {\Bo}} (Y) - \p_1 \Ncal_{ {\Bo}} (S_1) \nonumber \\
& = \p_1^2 \Ncal_\Bo(S_1) \left(y_1  - (- 1/2) \right) + \Rcal_*(Y)
\end{align*}
for all $Y= (y_1,y_2) \in B_{1/ 3 } (S_1) \setminus \GG_1$, where the second line holds thanks to \eqref{N_S_0_V_inequal}, and the third line does since $\p_1\p_2 \Ncal_\Bo(S_1) =0$.  The remainder $\Rcal_*(Y)$ satisfies
\beq\label{remainder}
\Rcal_*(Y)   \lesssim \left| Y-S_1 \right|^2.
\eeq
Since $\sin \Gb \Gt_1 > 0$ for $ \Gt_1 \in (0,\Ga) $, it follows from \eqref{a_e1-1}  that
\begin{align}
a_{1} &\le \frac{2}{(2\pi - \Ga){r_1^{\Gb}}} \Big[ \p_1^2  \Ncal_\Bo(S_1) \int_{0}^{2\pi-\Ga}  \left(y_1(r_1,\Gt_1)  - (- 1/2) \right) \sin \Gb_1 \Gt_1 d\Gt_1 \nonumber \\
&\qquad + \int_{0}^{2\pi-\Ga} \Rcal_*(Y) \sin \Gb_1 \Gt_1 d\Gt_1 \Big] \nonumber \\
&:= \frac{2}{(2\pi - \Ga)^{-1}  {r_1^{\Gb}}} \big[ \p_1^2 \Ncal_\Bo(S_1) I_1+ I_2 \big]. \label{last-1}
\end{align}
The above inequality holds for all $r_1 <1/3$. Thanks to \eqnref{remainder}, we have
$$
|I_2| \le C_2 r_1^2
$$
for some constant $C_2>0$. Since $\p_1^2 \Ncal_\Bo(S_1) < 0$ and
$$
I_1/r_1= - \int_{0}^{2\pi - \Ga} \cos \left( \Gt_1 + \frac {\Ga} 2 \right)  \sin \left(\Gb \Gt_1 \right) d\Gt_1 >0,
$$
we have
$$
\p_1^2 \Ncal_\Bo(S_1)  I_1 < - C_1 r_1
$$
for some constant $C_1>0$. By choosing $r_1$ sufficiently small, we infer from \eqnref{last-1} that $a_1 \neq 0$.

It then follows from \eqref{varphi-a_1B-temp} that
$$
|\nabla  \Gvf (Y) | \simeq |\nabla \Bcal_1 (Y)| \simeq  \frac 1 {|Y-S_1|^{1-\Gb}}
$$
for $Y \in B_{\Gk_0}(S_1) \setminus \overline{\GG_1}$ for some $\Gk_0$, from which \eqref{one-inclu} immediately follows. Now \eqref{one-inclu_main} is a simple consequence of \eqref{1st-ineq-one-inclusion} and \eqref{one-inclu}. \qed

\section{ Case 2: ${\Ba} = {\Bj } $ and $\Be = \Bo$}

In this section, we deal with the case when $\Ba$ and $\Be$, the direction and the location of the emitter, are given by
\beq\label{set-up-case2}
\Ba=\Bj=(0,1) \quad\mbox{and}\quad {\Be} = \Bo=(0,0).
\eeq
In this case, the equation in \eqref{def_gov_1} takes the form $\GD u = \p_2 \Gd_{\Bo}$. Since $\p_2 \Gd_{\Bo}$ satisfies, in the weak sense,
$$
\p_2 \Gd_{\Bo} (x_1,x_2) =  \p_2 \Gd_{\Bo} (-x_1,x_2) = -  \p_2 \Gd_{\Bo} (x_1,-x_2) ,
$$
we infer from the symmetry of the problem \eqref{def_gov_1} that
\beq\label{the_symmetries_Case_2}
u (x_1,x_2) = u (-x_1,x_2) =  -u(x_1,-x_2) .
\eeq
In particular, we have
\beq\label{equl_u_case2}
u|_{\p \GO_1}=  u|_{\p \GO_2}   = 0.
\eeq
Thus, in this case existence of corners on inclusions and the potential difference $u|_{\p \GO_2} -  u|_{\p \GO_1}$ do not contribute to the blow-up of $\nabla u$, only presence of emitter does. In fact, we prove the following theorem.

\begin{thm}\label{main2}
Let $u$ be the solution to \eqref{def_gov_1} under the assumption \eqref{set-up-case2}.  Then,
\beq\label{u_o_esti}
\left|\nabla  u (X)\right| \lesssim  \frac 1 {|X|^2}
\eeq
for all $X \in B_{1} \setminus  \left(\ol {\GO_{1} \cup \GO_{2}} \cup \{\Bo\} \right)$.
\end{thm}

We emphasize that
$$
\frac 1 {|X|^2} = 2\pi |\nabla \p_2 \Ncal_\Bo (X)|
$$
and
$$
\GD \p_2 \Ncal_\Bo = \p_2 \Gd_\Bo.
$$
So, \eqnref{u_o_esti} shows that $|\nabla u|$ is bounded by the field created by the emitter, and not enhanced by presence of the bow-tie structure.

By \eqref {def_c_a_e_e_2} and \eqref{equl_u_case2},  we have $ u = \Gs$. Thus we have from \eqnref{sigma_N_v} \beq\label{u_N+v_decomp}
u =  - v + \p_2 \Ncal_\Bo  \quad \mbox{in } \Rbb^2 \setminus (\ol{\GO_1 \cup \GO_1} \cup \{\Bo \}),
\eeq
where $v$ is the solution to \eqref{def_v}.

To prove Theorem \ref{main2}, we use \eqref{standard} (the standard gradient estimate). Thus we first obtain the following lemma for an upper bound of $v$.

\begin{lem}\label{v_0_lem}
It holds that
\beq
|v (X)| \leq \left| \p_2 \Ncal_\Bo (X)\right| \leq \frac 1 {2\pi} \frac 1 {|X|}
\eeq
for any $X \in \Rbb^2 \setminus {(\GO_1 \cup \GO_2 \cup \{\Bo\})}$.
\end{lem}
\pf
We follow the argument in the proof of Lemma \ref{lem_bounds}. Let $R_0 := (\Rbb  \times \Rbb^+) \setminus \overline{\GO_1 \cup \GO_2 }$ for ease of notation. Then we infer from \eqref{the_symmetries_Case_2} and \eqref{equl_u_case2} that
\beq\label{N-V_bdry_503}
\p_2 \Ncal_\Bo - v = 0   \quad\mbox{on } \p R_0 \setminus \{ \Bo \}.
\eeq
Since $\p_2 \Ncal_\Bo(X)=0$ if $x_2=0$, we have $v (X)=0$ for $X \in \p R_0 \cap \{ x_2=0 \}$. Moreover, $v > 0$ in $R_0$ since $\p_2 \Ncal_\Bo \ge 0$ on $\p R_0$. Thus, there are $r_0$ and $M$ such that
$$
Mx_2 >v(X) >0
$$
for any  $X \in \{ |X | < r_0 \} \cap R_0$. If $r < r_0$ is sufficiently small so that  $2\pi r^2  < M^{-1}$, then we have
\beq\label{NV_difference503}
\p_2 \Ncal_{\Bo}(X) = \frac{x_2}{2\pi |X|^2} > Mx_2 >v(X) >0
\eeq
for any  $X \in \{ |X | < r \} \cap R_0$. With \eqref {N-V_bdry_503} and \eqref{NV_difference503} in hand, the maximum principle yields
$$
\p_2 \Ncal_{\Bo} - v >0  \quad\mbox{in } R_0 \cap \{ |X| \ge r \}.
$$
Since $r$ can be arbitrarily small, we conclude that
$$
\p_2 \Ncal_{\Bo} > v >0  \quad\mbox{in } R_0 .
$$

Similarly, one can show that
$$
\p_2 \Ncal_{\Bo } < v <0 \quad\mbox{in }  (\Rbb  \times \Rbb^-) \setminus \ol{(\GO_1 \cup \GO_2 )} ,
$$
and the proof is complete.
\qed

\medskip
We are now ready to prove Theorem \ref {main2}.

\noindent{\sl Proof of Theorem \ref {main2}}. Let $R:= B_{1}\setminus \left( \ol{\GO_1 \cup \GO_2} \cup \{\Bo\} \right)$ and decompose $R$ into four disjoint regions:
\begin{align*}
R_1 &:= B_{3\Ge  /8 } \setminus \{ \Bo \}, \\
R_2 &:= \left(B_{\Ge / 4 } (V_1) \cup B_{\Ge / 4 } (V_2) \right) \cap R, \\
R_3 &:= (B_{3\Ge} \cap R )\setminus \left( R_1 \cup R_2 \right) , \\
R_4 &:=  R \setminus (R_1 \cup R_2 \cup R_3) = R \setminus  B_{3\Ge}.
\end{align*}
The desired estimates for $\nabla u$ in regions $R_1$, $R_3$ and $R_4$ are derived by \eqnref{standard} using Lemma \ref{v_0_lem}. However, the estimate in $R_2$ is somewhat different.

We first deal with the estimates in $R_1$, $R_3$ and $R_4$.
Every point $X$ in $R_1$ satisfies $ B_{|X|/3} (X) \subset R$. By \eqref  {u_N+v_decomp} and Lemma \ref{v_0_lem},  \eqnref{standard} yields \begin{align}  |\nabla u(X)| & \leq | \nabla v(X) |  + |\nabla \p_{2} \Ncal_\Bo (X)|  \notag\\& \leq \frac  2 {|X|/3}  {\norm {v}_{L^{\infty}(B_{|X|/3} (X)) }} +  |\nabla \p_{2} \Ncal_\Bo (X)| \notag \\& \leq \frac  6 {|X|}  {\norm {\p_{2} \Ncal_\Bo}_{L^{\infty}(B_{|X|/3} (X)) }} +  |\nabla \p_{2} \Ncal_\Bo (X)| \lesssim  \frac 1 {|X|^2}. \label{R_1_gradi_estimate1}
\end{align}

Note that for any point $X \in R_3$, the distance between $X$ and $\{V_1, V_2, \Bo\}$ is greater than $\Ge/4$. Since $u = 0$ on $\p \GO_1 \cup \p \GO_2$, the function $u$ can be extended by reflection into $B_{\Ge / 8 } (X)$ as a harmonic function. If we denote the extended function by $u$, then Lemma \ref{v_0_lem} implies that
\beq\label{R_3_stnd_grd_esti}
\norm {u}_{L^{\infty} (B_{\Ge/8} (X))} = \norm {u}_{L^{\infty} (B_{\Ge/8} (X) \cap R)}  \leq  \norm {\p_2 \Ncal_\Bo}_{L^{\infty} (B_{\Ge/8} (X))} .
\eeq
It follows from \eqnref{standard} that
\beq\label{R_3_stnd_grd_esti2}
 |\nabla u (X)|  \leq  \frac 2  {\Ge/8} \norm {u}_{L^{\infty} (B_{\Ge/8} (X))}\leq  \frac {16}  {\Ge} \norm {\p_2 \Ncal_\Bo}_{L^{\infty} (B_{\Ge/8} (X))}  \lesssim \frac 1 {|X|^2},
\eeq
where the last inequality holds to be true since $|X|< 3\Ge$.

To deal with the estimate in $R_4$, let $\mu$ be the constant appeared in \eqref{GOone-2} and choose a positive constant $c < 1/3$ such that $c |X| <  \mu - 1$ for all $X \in R_4$. Then for $X \in R_4$, we have
$$
B_{c |X|} (X)\subseteq B_{\mu}\setminus B_{2\Ge}.
$$
Thus, the function $u$ can be extended into $B_{c |X| } (X)$ as a harmonic function, which is still denoted by $u$, since $u = 0 $ on $\p \GO_1 \cup \p \GO_2$. In the same way as before, Lemma \ref{v_0_lem} yields
$$
\norm {u}_{L^{\infty} (B_{c |X|} (X))} = \norm {u}_{L^{\infty} (B_{c |X|} (X) \cap R)}  \leq  \norm {\p_2 \Ncal_\Bo}_{L^{\infty} (B_{c |X|} (X))} \leq \frac 1 {2\pi(1-c)} \frac 1 {|X|},
$$
and \eqnref{standard} yields
$$
|\nabla u (X)|  \leq  \frac 2  {c|X|} \norm {u}_{L^{\infty} (B_{c|X|} (X))} \leq \frac 1 {\pi c (1-c)} \frac 1  {|X|^2}.
$$

To estimate $\nabla v$ in $R_2$, we first observe that $R_2 = (B_{\Ge/4}(V_1) \setminus \GO_1) \cup (B_{\Ge/4}(V_2) \setminus \GO_2)$. Since $\p B_{3\Ge /8 } (V_j) \setminus \GO_j \subseteq R_1 \cup R_3$ ($j=1,2$), we infer from  \eqref{R_1_gradi_estimate1} and  \eqref{R_3_stnd_grd_esti2} that
\beq\label{linfty}
\norm {\nabla u}_{L^{\infty} \left( \p B_{3\Ge /8 } (V_1) \setminus \GO_1  \right)} \lesssim \frac 1 {\Ge^2 } .
\eeq
Thanks to \eqnref{the_symmetries_Case_2} and \eqnref{equl_u_case2}, $u$ admits the following Fourier series expansion on $X \in \ol{B_{3\Ge /8} (V_j)} \setminus \GO_j$:
\beq
u(X) = \sum_{ n=1 } ^{\infty} {a}_{2n} |X-V_j|^{2n\Gb} \sin\left(2n \Gb \Gt_j ( X/\Ge) \right) ,
\eeq
where $\Gt_j$ is defined in \eqref{angle}. By taking the radial and angular derivatives, one can see that
\beq\label{gradpolar}
\left| \nabla u (X)\right| \lesssim \sum_{n=1}^{\infty} {2n \Gb } \left|  {a}_{2n}\right|  |X-V_j| ^{2n \Gb -1}.
\eeq
We then obtain from the Cauchy-Schwartz inequality that if $X \in B_{\Ge/4}(V_j) \setminus \GO_j$, then
\begin{align*}
\left| \nabla u (X)\right| &\leq \left( \sum_{n=1}^{\infty} {4n^2 \Gb^2 }   |{a}_{2n}|^2  \left(\frac {3\Ge}{8} \right)^{4n \Gb -2}\right)^{\frac 1 2 } \left(\sum_{n=1} \frac {|X-V_j|^{4n \Gb -2}} {\left(3\Ge / 8 \right)^{4n \Gb -2}}\right)^{\frac 1 2} \\
&\lesssim \left( \sum_{n=1}^{\infty} {4n^2 \Gb^2 } |{a}_{2n}|^2  \left(\frac {3\Ge}{8} \right)^{4n \Gb -2}\right)^{\frac 1 2 },
\end{align*}
where the last inequality holds to be true since $4n\Gb -2 >  0 $ for $n=1,2,3, \ldots$. We then have from \eqnref{linfty} and \eqnref{gradpolar} that
\begin{align*}
\sum_{n=1}^{\infty} {4n^2 \Gb ^2 } |{a}_{2n}|^2  \left(\frac {3}{8} \Ge \right)^{4n \Gb -2}
& \lesssim \frac{1}{\Ge} \norm {\nabla u}_{L^{2} \left( \p B_{3\Ge /8 } (V_j)\setminus \GO_j \right)}^2 \\
& \lesssim \norm {\nabla u}^2 _{L^{\infty} \left( \p B_{3\Ge /8 } (V_j)\setminus \GO_j \right)} \lesssim  \frac 1 {\Ge ^4 } .
\end{align*}
Since $\Ge \simeq |X|$ for all $X \in B_{ \Ge /4 } (V_j) \setminus \GO_j$, we have
$$
\left| \nabla u (X)\right| \lesssim  \frac 1 {|X|^2}.
$$
This completes the proof of Theorem \ref{main2}. \qed

\section{Case 3: ${\Ba} = {\Bj } $ and $\Be \neq \Bo$ }

In this section, we deal with the case when
\beq\label{set-up-case3}
\Ba =\Bj \quad \mbox{and} \quad \Ge\Be \neq \Bo,
\eeq
so that the equation in \eqref{def_gov_1} takes the form $\GD u = \p_2 \Gd_{\Ge\Be}$. So, similarly to \eqnref{the_symmetries_Case_2}, we see that
\beq
u (x_1,x_2) =  u (-x_1,x_2),
\eeq
which implies, in particular, that there is no potential difference, namely,
\beq\label{nodiff}
u |_{\p \GO_2} - u |_{\p \GO_1} = 0.
\eeq
Thus the blow-up of $\nabla u$ is only caused by presence of corners and the emitter.

We establish two upper bounds of $\nabla u$, one in a neighborhood of vertices of the inclusions and the other in the rest (Theorem \ref{main3-1}). We also derive a lower bound in a neighborhood of vertices of the inclusions (Theorem \ref{main3-2}), under a mild geometric assumption on inclusions (see the condition (A) below).

\begin{thm}\label{main3-1}
Let $u$ be the solution to \eqref{def_gov_1} under the assumption \eqref{set-up-case3}, and
\beq\label{Rdef}
R := B_1 \setminus \overline {\GO_1 \cup \GO_2 \cup \{\Ge {\bf e}\}}.
\eeq
There exists a positive constant $\Ge_0$ such that for a given $c_0 \in (0,1/2)$ the following estimates hold with a constant $C$ for all $\Ge < \Ge_0$:
\begin{itemize}
\item[(i)]   For all $X \in  \left( B_{c_0\Ge}(V_1) \cup  B_{c_0\Ge}(V_2)  \right) \cap R$,
\beq\label{upper_bd_main3_upp1}
|\nabla u (X)| \le \frac C {\Ge^{1+\Gb}} \sum_{j=1} ^2 \frac 1 {|X- V_j|^{1-\Gb}} .
\eeq
\item[(ii)] For all $X \in R \setminus ( B_{c_0\Ge}(V_1) \cup B_{c_0\Ge}(V_2)) $,
\beq\label{upper_bd_main3_upp2}
\left|\nabla u (X)\right|\le \frac C {|X - \Ge {\bf e}|^2} .
\eeq
\end{itemize}
\end{thm}
Note that $|X - \Ge {\bf e}|^{-2} \simeq \left| \nabla  \p_{2} \Ncal_{\Ge {\bf e}} (X) \right|$. Therefore, \eqnref{upper_bd_main3_upp2} shows that the blow-up of $|\nabla u|$ away from the vertices is caused by presence of the emitter, not by presence of corners. On the other hand, \eqnref{upper_bd_main3_upp1}, which is similar to \eqnref{utwonear}, shows that $\nabla u$ may be enhanced because of presence of corners.

\begin{figure}[h!]
\begin{center}
\epsfig{figure=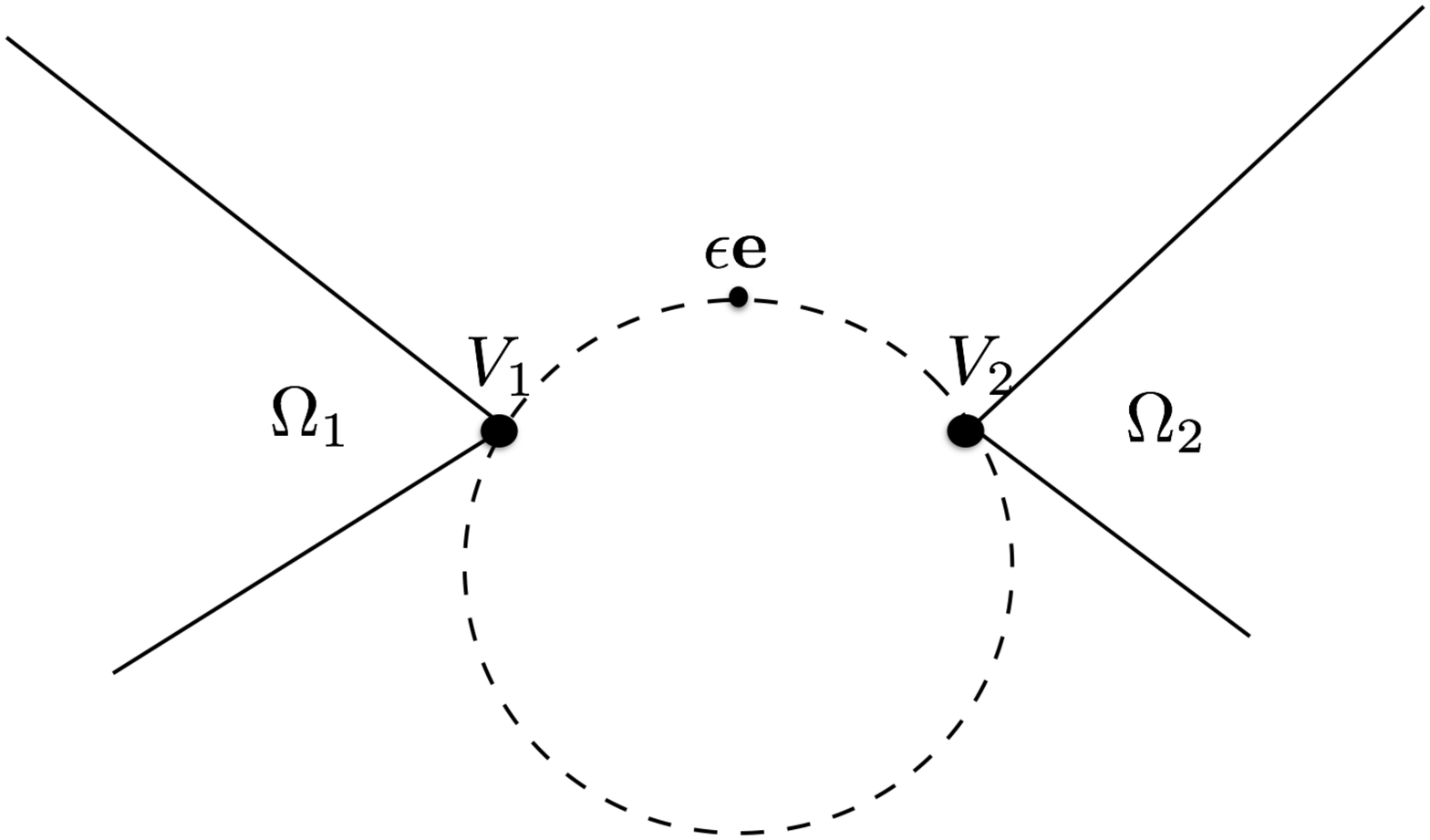,width=7cm}
\end{center}\vskip -20pt
\caption{The condition (A) on $\GO_1$, $\GO_2$ and $\Ge {\Be}$}\label{Fig1}
\end{figure}

The following theorem shows that field enhancement actually occurs under a mild assumption on $\Be$, $\GG_1$, and $\GG_2$:
\begin{itemize}
\item[(A)] the circle passing through three points $S_1$, $S_2$ and $\Be$ does not meet with $\p \GG_1  \setminus \{S_1\} $ and $\p \GG_2 \setminus \{S_2\}$.
\end{itemize}
After the scaling we may rephrase the condition (A) as follows: the circle passing through three points $V_1$, $V_2$ and $\Ge {\Be}$ does not meet with $\p \GO_1  \setminus \{V_1\} $ and $\p \GO_2 \setminus \{V_2\}$ (see Figure \ref{Fig1}). It is worthwhile to mention that it may fail to hold, for example, if the aperture angle of the vertices is large and $p \neq \pm 1/2$.

\begin{thm}\label{main3-2}
Let $u$ be the solution to \eqref{def_gov_1} under the assumption \eqref{set-up-case3}.
If the condition (A) is fulfilled, then there exist  positive constants $\tilde c_0$ and $\Ge_0$ such that
\beq\label{lower_bd_main3}
|\nabla u (X)| \gtrsim  \frac 1 {\Ge^{1+\Gb}} \sum_{i=1} ^2 \frac 1 {|X- V_i|^{1-\Gb}}
\eeq  for all $X \in  \left( B_{\tilde c_0 \Ge  }(V_1) \cup  B_{\tilde c_0 \Ge  }(V_2)  \right) \cap R$ ($R$ is defined in \eqnref{Rdef}) and $\Ge < \Ge_0$.
\end{thm}

Proofs of these theorems rely on the following three lemmas whose proofs are given in subsection \ref{subsec:proof3}.

\begin{lem}\label{v_V_boundedness}
Let  $v$ and $V$ be  the solutions to \eqref {def_v} and  \eqref{def_V}, respectively. There exists a positive constant $ \Ge_0$ such that the following holds for all $\Ge< \Ge_0$
\beq\label{vCase3}
|v (X)|  \lesssim  \frac 1 { \Ge + |X|}
\eeq
for all $X  \in \Rbb^2 \setminus \ol {(\GO_{1} \cup \GO_{2})}$, and
\beq\label{VCase3}
\frac 1 {\Ge}\left |V\left( {\Ge}^{-1}Y\right)\right|   \lesssim   \frac 1 {\Ge + |Y|}
\eeq
for all $Y \in \Rbb^2 \setminus \ol {(\GG_{1} \cup \GG_{2})}$.
\end{lem}

\begin{lem}\label{thm23_lema1}
Let $\Gvf$ be the function defined in \eqref{defGvf}. There exists a constant $a_{1}$ such that for every positive $\Gk <\frac 1 2$ there is a constant $C$ satisfying
\beq\label{case3_variphi_B_1_B_2}
\left| \nabla \Gvf (Y) - a_{1} \left ( \nabla \Bcal_1 (Y) +  \nabla \Bcal_2 (Y) \right) \right | \leq C
\eeq
for all $Y \in B_{\Gk} (S_1) \cup B_{\Gk} (S_2)$. Moreover, if the condition (A) is fulfilled, then $a_1 \neq 0$.
\end{lem}

\begin{lem} \label{lem:ACC}
Let $\Rcal_2$ be as defined in \eqref{Rtwo}. There exists a positive constant $\Ge_0$ such that
\beq \left| \nabla  \Rcal_2(X)\right| \lesssim  \frac 1 {\Ge^{1-\Gb}} \left(  \left|(\nabla \Bcal_1)\left( {\Ge}^{-1} X\right) \right| + \left| (\nabla \Bcal_2)\left( {\Ge}^{-1} X\right)\right| \right)  \label{lem:ACC_ineq} \eeq for all $X \in B_{1  } \setminus \ol {(\GO_{1}\cup \GO_{2})}$ and for all $\Ge < \Ge_0$.
\end{lem}

\noindent{\sl Proofs of Theorems \ref{main3-1} and \ref{main3-2}}.
In view of \eqref{uone}, \eqref{Rone} and \eqnref{nodiff},  we have $u_1 = \Rcal_1 \equiv 0$. Thus we see from \eqref{utwo} and \eqref{udecom} that
\beq
\nabla u (X) = \nabla u_2 (X) + \nabla  \Rcal_2 (X)= \Ge^{-2}(\nabla \Gvf) (\Ge^{-1}X) + \nabla  \Rcal_2 (X)
\eeq
for any $X \in R$, where $R$ is the region defined by \eqnref{Rdef}.

Let $\Ge_1$ be the smallest of $\Ge_0$ appearing in Lemmas \ref {v_V_boundedness} and \ref{lem:ACC}. Then, \eqnref{upper_bd_main3_upp1} results from  \eqref {case3_variphi_B_1_B_2} and \eqref{lem:ACC_ineq} as  follows:
\begin{align*}
|\nabla u (X)|& \lesssim  \left| a_{1} {\Ge^{-2}} \left ( ( \nabla \Bcal_1 )\left({\Ge }^{-1}X\right)   +( \nabla \Bcal_2 )\left({\Ge }^{-1}X\right)     \right)    \right| + \Ge^{-2} + |\nabla \Rcal_2 (X) | \\
& \lesssim  {\Ge^{-2}}  \left|  ( \nabla \Bcal_1 )\left({\Ge }^{-1}X\right) \right|   + {\Ge^{-2}}   \left| ( \nabla \Bcal_2 )\left({\Ge }^{-1}X\right)   \right| \\
& \lesssim  \frac 1 {\Ge^{1+\Gb}} \sum_{j=1} ^2 \frac 1 {|X- V_j|^{1-\Gb}}
\end{align*}
for all $X \in  \left( B_{c_0\Ge}(V_1) \cup  B_{c_0\Ge}(V_2)  \right) \cap R$ and all $\Ge < \Ge_1$, where the last inequality follows from \eqnref{nablaBcal}.

If the condition (A) is satisfied, then we infer from Lemmas \ref{thm23_lema1} and \ref{lem:ACC} that there are three positive numbers $\Ge_0 < \Ge_1$,  $\tilde c_0$ and $C$ such that
\begin{align*}
|\nabla u (X)|& \gtrsim  \left| a_{1} {\Ge^{-2}} ( \nabla \Bcal_1 )\left({\Ge }^{-1}X\right) \right|   -\left| a_{1} {\Ge^{-2}}  ( \nabla \Bcal_2 )\left({\Ge }^{-1}X\right) \right|  - C \left( \Ge^{-2} + |\nabla \Rcal_2 (X) | \right) \\
& \gtrsim  \left(|a_1| \frac 1 {\Ge^{1+\Gb}}  - C \right)\frac 1 {|X-V_1|^{1-\Gb}} - \left(|a_1| \frac 1 {\Ge^{1+\Gb}}  + C \right)\frac 1 {|X-V_2|^{1-\Gb}} - C \Ge^{-2} \\
& \gtrsim  \frac 1 {\Ge^{1+\Gb}}  \frac 1 {|X-V_1|^{1-\Gb}} \gtrsim \frac 1 {\Ge^{1+\Gb}} \sum_{j=1} ^2 \frac 1 {|X- V_j|^{1-\Gb}}
\end{align*} for all $X \in  B_{\tilde c_0 \Ge  }(V_1) \cap R$ and all $\Ge < \Ge_0$. We obtain the same inequality for $X \in  B_{\tilde c_0 \Ge  }(V_2) \cap R$, and hence, \eqref{lower_bd_main3} follows.

We now prove \eqref{upper_bd_main3_upp2}. Since $u_1 = \Rcal_1 = 0$, we see from \eqref{def_v} that
\beq\label{u-one}
u -u|_{\p \GO_1} =\Gs -\Gs|_{\p \GO_1} =  \p_{2} \Ncal_{\Ge \bf e } - v \quad \mbox{in } R .
\eeq
Let
\beq\label{Ronetwo}
R_1:= R \setminus \overline{B_{2\Ge} } \quad\mbox{and}\quad R_2:=(B_{2\Ge} \cap R )\setminus ( B_{ c_0 \Ge  }(V_1) \cup B_{c_0 \Ge  }(V_2)),
\eeq
so that
$$
R \setminus ( B_{ c_0 \Ge  }(V_1) \cup B_{ c_0 \Ge   }(V_2)) = R_1 \cup R_2.
$$

Let $c_1 = \min (\mu -1, 1/2)$, where $\mu$ is given in \eqref{GOone}. Then one can see that for any $X \in R_1$, $B_{c_1 |X|} (X) \subset B_{\mu} \setminus B_{\Ge}$ which  does not contain any of $\Ge \bf e $, $V_1$ and $V_2$. Since $u - u |_{\p \GO_1} = 0$ on $\p \GO_1\cup \p \GO_2$, $u - u |_{\p \GO_1}$ can be extended by reflection to $B_{c_1 |X|} (X)$ as a harmonic function. Denote the extended function by $U$. Then we have
\beq
\norm{U}_{L^{\infty} (B_{c_1 |X|} (X))} \le 2 \norm {u - u |_{\p \GO_1} }_{L^{\infty} (B_{c_1 |X|} (X)\cap R)}.
\eeq
By \eqnref{standard} we have
$$
|\nabla u (X)| =|\nabla U(X)| \le \frac {2}{c_1 |X| } \norm{U}_{L^{\infty}(B_{c_1 |X|} (X))} \le \frac {4}{c_1 |X| } \norm {u - u |_{\p \GO_1} }_{L^{\infty} (B_{c_1 |X|} (X)\cap R)}.
$$
It then follows from \eqnref{vCase3} and \eqnref{u-one} that
\begin{align}
|\nabla u (X)| &  \leq \frac {4}{c_1 |X| } \norm { \p_2 \Ncal_{\Ge \bf e} }_{L^{\infty} (B_{c_1 |X|} (X)\cap R)} + \frac {4}{c_1 |X| } \norm { v }_{L^{\infty} (B_{c_1 |X|} (X)\cap R)} \notag \\&\lesssim \frac 1 {|X|^2} \lesssim \frac 1 {|X - \Ge \bf e|^2}  \label {case3_last_ineq_former_region}
\end{align}
for any $\Ge < \Ge_0$, since $|{\bf e}| < 1$, if $X \in R_1$.

To deal with the second region $R_2$, we first observe from the mean value theorem that if $X \in (\p B_{3\Ge}  ) \cap R$, then there is $X^*\in (\p B_{3\Ge}  ) \cap R$ such that
$$
|u(X) - u |_{\p \GO_1}| \lesssim \Ge |\p_t u(X^*)|,
$$
where $\p_t$ denotes the tangential derivative along $\p B_{3\Ge}$. Thus we have
$$
\norm{u - u |_{\p \GO_1}}_{L^{\infty}((\p B_{3\Ge})\cap R)} \lesssim \Ge\norm{\nabla u }_{L^{\infty}((\p B_{3\Ge})\cap R)}  \lesssim \Ge^{-1} ,
$$
where the last inequality follows from \eqref {case3_last_ineq_former_region} since $(\p B_{3\Ge})\cap R \subset R_1$.
It then follows from \eqnref{u-one} that
$$
\norm{v}_{L^{\infty}( (\p B_{3\Ge}) \cap R )}\leq \norm{\p_2 \Ncal_{\Ge {\bf e}} }_{L^{\infty}( (\p B_{3\Ge}) \cap R )} + \norm{u - u |_{\p \GO_1}}_{L^{\infty}(  (\p B_{3\Ge}  ) \cap R )} \lesssim \Ge^{-1} .
$$
By the boundary condition of $v$, we have
$$
\norm{v}_{L^{\infty}( ( \p R ) \cap B_{3\Ge})}= \norm{\p_ 2\Ncal_{\Ge {\bf e}}}_{L^{\infty}( ( \p R ) \cap B_{3\Ge} )} \lesssim \Ge^{-1}.
$$
Note that $v$ is defined in $B_{3\Ge} \cap (R \cup \{\Ge {\bf e }\} )$. So, we infer from the maximum principle that
\beq\label{case3_last_ineq_latter_region}
\norm{v}_{L^{\infty}(  B_{3\Ge} \cap R )} \lesssim \Ge^{-1}.
\eeq

On the other hand, there is a positive constant $c_2 < 1/2$ independent of $\Ge$ such that the following holds for every $X \in R_2$:
$$
B_{c_2|X-\Ge {\bf e} |} (X) \subset B_{3\Ge} \setminus ( B_{ c_0 \Ge /2  }(V_1) \cup B_{c_0  \Ge /2 }(V_2) \cup \{\Ge {\bf e}\}  ).
$$
As before, $u - u |_{\p \GO_1}$ can be extended to $B_{c_2|X-\Ge {\bf e} |} (X) $ as a harmonic function, and the extended function, which we denote by $U_1$, satisfies
$$
\norm{U_1}_{L^{\infty} (B_{c_2|X-\Ge {\bf e} |} (X))} \le 2 \norm{ u - u |_{\p \GO_1}}_{L^{\infty} (B_{c_2|X-\Ge {\bf e} |} (X) \cap R)}.
$$
Then \eqnref{standard} and \eqref{case3_last_ineq_latter_region} yield that
\begin{align*}
|\nabla u (X)| &=|\nabla \left(u (X) - u |_{\p \GO_1}\right)| \leq \frac {4}{c_2 |X-{\Ge {\bf e}}| } \norm {u - u |_{\p \GO_1} }_{L^{\infty} (B_{c_2 |X-{\Ge {\bf e}}|  } (X)\cap R)} \notag\\&  \leq \frac {4}{c_2 |X-{\Ge {\bf e}}|  } \norm { \p_2 \Ncal_{\Ge \bf e} }_{L^{\infty} (B_{c_2 |X-{\Ge {\bf e}}|  } (X)\cap R)} + \frac {4}{c_2 |X-{\Ge {\bf e}}|  } \norm { v }_{L^{\infty} (B_{c_2 |X-{\Ge {\bf e}}|  } (X)\cap R)} \notag \\&\lesssim \frac 1 {|X- \Ge {\bf e}|^2} +  \frac 1 {\Ge |X- \Ge {\bf e}| } \lesssim \frac 1 {|X - \Ge \bf e|^2}  \notag
\end{align*}
for all $\Ge < \Ge_0$. Therefore, we obtain the \eqref{upper_bd_main3_upp2} in $R_2$ and the proof is complete.
\qed

\subsection {Proofs of Lemma \ref{v_V_boundedness}, \ref{thm23_lema1} and \ref{lem:ACC}}\label{subsec:proof3}

\begin{figure}[h!]
\begin{center}
\epsfig{figure=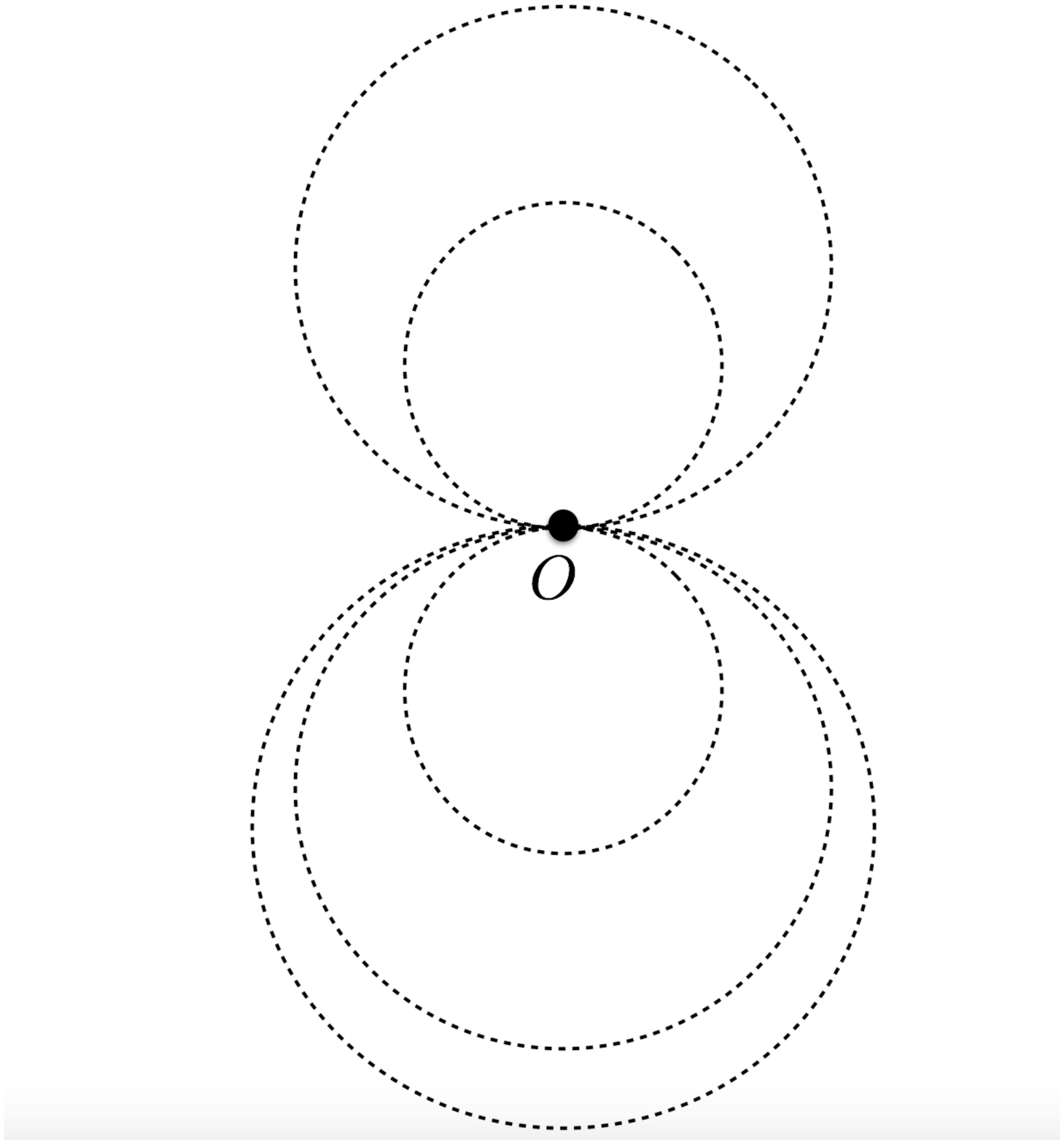,width=5cm}
\end{center} \vskip -30pt
\caption{Level curves of $\p_2 \Ncal_{\Bo}$}\label{Fig2}
\end{figure}

One can easily see that for any real number $\Gz \neq 0$,
$$
\p_2 \Ncal_\Bo (x_1,x_2) = \frac {x_2}{2\pi (x_1^2 + x_2^2)} = \frac 1 {4\pi \Gz}
$$
if and only if
$$
x_1^2 + (x_2-\Gz)^2 = \Gz^2 \mbox{ and } (x_1,x_2)\neq \Bo.
$$
We can infer from this that the level curves of $\p_2 \Ncal_{\Bo}$ are circles passing through $\Bo$ with the center on the $y$-axis, as depicted in Figure \ref{Fig2}. The circular level set with the center above (below, resp.) $\Bo$ corresponds to the positive (negative, resp.) value. The larger the radius is, the smaller is the value in absolute.
Since $\p_2 \Ncal_{\Ge\Be}$ is nothing but a translate of $\p_2 \Ncal_\Bo$, we have the following lemma.

\begin{lem}\label{lem_circle_passing}
The level curves of $\p_2 \Ncal_{\Ge\Be}$, which is defined on $\Rbb^2 \setminus \{\Ge \Be \}$, are circles passing through $\Ge\Be$ with the center on the $y$-axis. The circular level set with the center above (below, resp.) $\Ge\Be$ corresponds to the positive (negative, resp.) value. The larger the radius is, the smaller is the value in absolute.
\end{lem}

Recall that $\Ge {\Be} = \Ge (0,p)$. We assume $p>0$ without loss of generality.

\begin{lem}\label{lem_mini_v_E}
Let  $v$ and $V$ be  the solutions to \eqref {def_v} and  \eqref{def_V}, respectively.   There exists a positive constant $\Ge_0$ such that
\begin{align}
& \max \left\{ v (X) ~|~ X \in \Rbb^2 \setminus {(\GO_{1} \cup \GO_{2})}  \right\} \simeq  \Ge^{-1}, \label{vone} \\
& \min \left\{ v (X)~ |~ X \in \Rbb^2 \setminus {(\GO_{1} \cup \GO_{2})}  \right\} \simeq  - \Ge^{-1}, \label{vtwo}
\end{align}
for all $\Ge < \Ge_0$, and
\begin{align}
& \max \left\{ V (Y) ~|~ Y \in \Rbb^2 \setminus {(\GG_{1} \cup \GG_{2})}  \right\} \simeq 1, \label{Vone}\\
& \min \left\{ V (Y) ~|~ Y \in \Rbb^2 \setminus {(\GG_{1} \cup \GG_{2})}  \right\} \simeq -1. \label{Vtwo}
\end{align}
Moreover, if the condition (A) is fulfilled, then the function $v$ has the negative minimum value $\p_2 \Ncal_{\Ge {\Be}} (V_1) = \p_2 \Ncal_{\Ge {\Be}} (V_2)$  at $V_1$ and $V_2$ for $\Ge< \Ge_0$, and similary the function $V$ has the negative minimal value $\p_2 \Ncal_{ {\Be}} (S_1) = \p_2 \Ncal_{{\Be}} (S_2)$  at $S_1$ and $S_2$.
\end{lem}

\pf
There are two circles passing through $\Be$, one of whose centers is on the $y$-axis and above $\Be$, and the other is on the $y$-axis and below $\Be$, such that their radii are the smallest under the condition that circles meet $\p\GG_1$ and $\p\GG_2$ above $\Be$ and below $\Be$, respectively. Let $C_a$ be the circle above $\Be$ and $C_b$ be the circle below $\Be$.  Then $C_a$ osculates the $\p\GG_1$ and $\p\GG_2$ in the upper plane, say at $P_1$ and $P_2$, respectively. The circle $C_b$ meets with $\p\GG_1$ and $\p\GG_2$ at say $Q_1$ and $Q_2$, respectively.

\begin{figure}[h!]
\begin{center}
\epsfig{figure=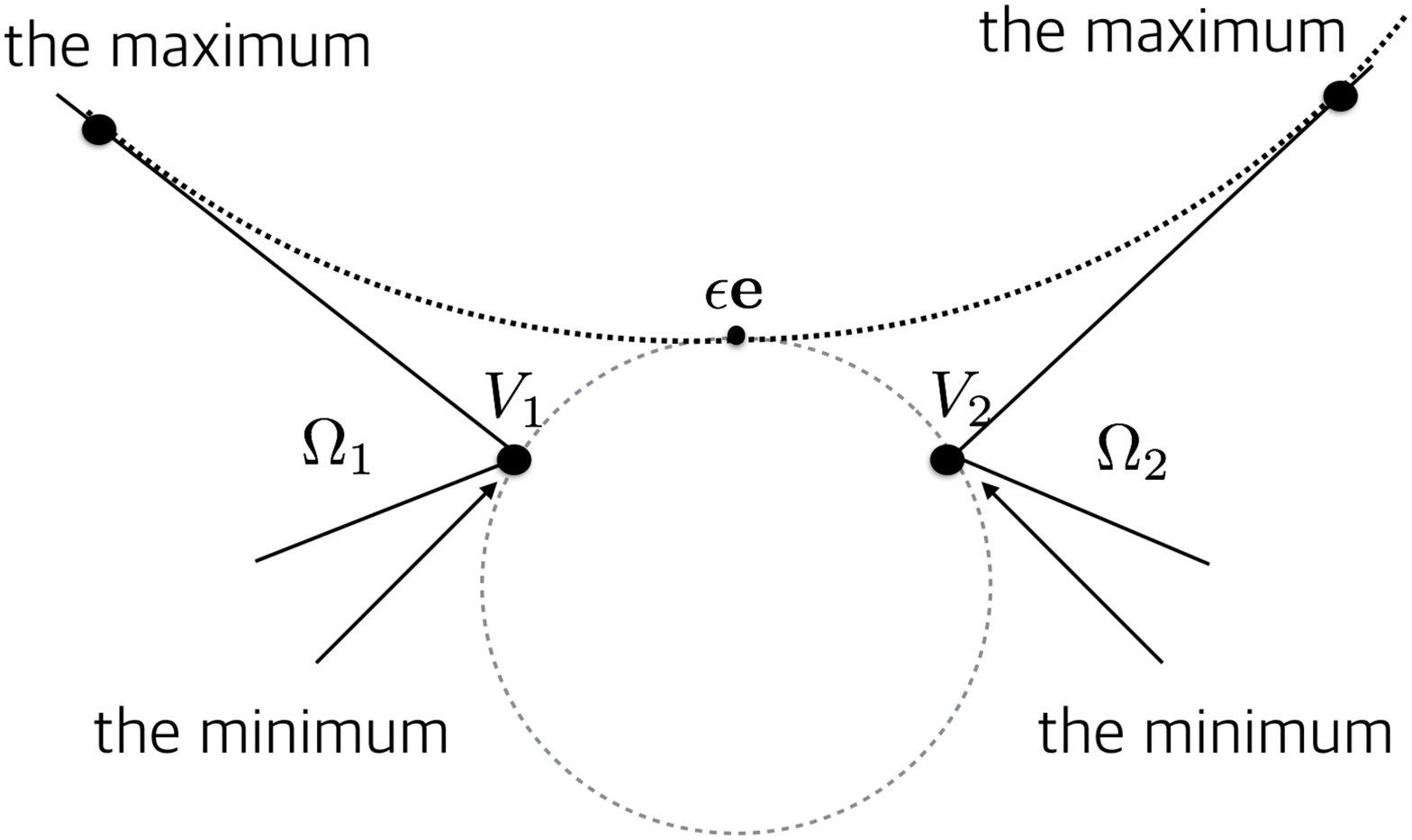,width=7cm}
\end{center} \vskip -20pt
\caption{Level curves and the maximum and minimum of $\p_{2} \Ncal_{\Ge {\Be}}\big|_{\p \GO_1 \cup \p \GO_2}$.}\label{Fig3}
\end{figure}

The scaled circles $\Ge C_a$ and $\Ge C_b$ meet with $\p \GO_1$ and $\p\GO_2$ at $\Ge P_1 \in \p \GO_1 $ and $\Ge P_2 \in \p \GO_2$, and $\Ge Q_1 \in \p \GO_1 $ and $\Ge Q_2 \in \p \GO_2$, respectively, provided that $\Ge$ is so small that all of these four points lie in the unit disk. See Figure \ref{Fig3}.
The circle $\Ge C_a$ is the smallest circle above $\Ge\Be$ passing through $\Ge\Be$. So, by Lemma \ref{lem_circle_passing}, we infer that $\p_{2} \Ncal_{\Ge {\Be}}$ restricted to $\p \GO_1 \cup \p \GO_2$ attains its maximum value at $\Ge P_1$ and $\Ge P_2$. Similarly, we see that it attains its minimum value  at $\Ge Q_1$ and $\Ge Q_2$.
Since $v=\p_{2} \Ncal_{\Ge {\Be}}$ on $\p \GO_1 \cup \p \GO_2$, $v$ as a function in $\Rbb^2 \setminus (\GO_1 \cup \GO_2)$ attains the same maximum and minimum values at the same points. By \eqref{scaling}, the maximum value is $ {\Ge}^{-1} \p_2 \Ncal_{\Be} (P_1) =  {\Ge}^{-1}\p_{2} \Ncal_{\Be} (P_2)$, and the minimum value is $ {\Ge}^{-1} \p_2 \Ncal_{\Be} (Q_1) =  {\Ge}^{-1} \p_2 \Ncal_{\Be} (Q_2)$. Since $\p_2 \Ncal_{\Be} (P_1) >0$ and $\p_2 \Ncal_{\Be} (Q_1) <0$ regardless of $\Ge$, \eqnref{vone} and \eqnref{vtwo} hold. One can prove \eqnref{Vone} and \eqnref{Vtwo} similarly.

If the condition (A) is fulfilled, then $\Ge C_b$ is the circle passing through three points $V_1$, $V_2$ and $\Ge {\Be}$. So, in this case, $Q_j=S_j$ and $\Ge Q_j=V_j$ for $j=1,2$. Thus we have the last sentence in the lemma.  \qed

\medskip
\noindent{\sl  Proof of Lemma \ref{v_V_boundedness}}.
Let
$$
H_+:= \left(\Rbb \times (0,\infty)\right)\setminus \ol {(\GO_{1} \cup \GO_{2})} \quad\mbox{and}\quad H_-:= \left(\Rbb \times (-\infty,0)\right)\setminus \ol {(\GO_{1} \cup \GO_{2})}.
$$
We still assume that $p>0$, namely, $\Ge\Be \in H_+$.

There is a constant $C$ depending on $p$ such that
$$
| v (X) | = | \p_2 \Ncal_{\Ge {\Be}} (X)| \leq  C \p_2 \Ncal_{-\Ge {\Be}}(X)
$$
for  $ X \in\p H_+$.
Since $\p_2 \Ncal_{-\Ge {\Be}}$ is a positive harmonic function in $H_+$, we infer from the maximum principle
$$
| v (X) |  \leq  C \p_2 \Ncal_{-\Ge {\Be}} (X) \leq \frac{C}{\Ge + |X|}
$$
for  $X \in  \overline{H_+}$. Here and in the following, $C$ stands for constants depending on $p$ which may differ at each occurrence.

We also have
$$
| v (X) | = | \p_2 \Ncal_{\Ge {\Be}} (X)| \leq  -\p_2 \Ncal_{\Ge {\Be}}(X)
$$
for  $ X \in\p H_-$.
Since $-\p_2 \Ncal_{\Ge {\Be}}$ is a positive harmonic function in $H_-$, we have
$$
| v (X) |  \leq  -\p_2 \Ncal_{\Ge {\Be}} (X) \leq \frac{C}{\Ge + |X|}
$$
for  $X \in \overline{ H_-}$. Thus \eqnref{vCase3} follows.

The estimate \eqnref{VCase3} can be proved similarly.
\qed

\medskip

\noindent{\sl Proof of Lemma \ref{thm23_lema1}}.
Let $\Gt_1$ and $r_1$ be the polar coordinates with respect to $S_1$ as defined in \eqref{angle} and \eqref{rldef}. Since $\Gvf$ vanishes  on $\p \GG_1$, it admits the Fourier series expansion of the form
\beq\label{6000}
\Gvf (Y) = \sum_{n=1}^{\infty} a_{n} r_1^{n\Gb} (Y) \sin n\Gb \Gt_1 (Y)
\eeq
for $ Y \in B_{ 1/2} (S_1) \setminus \ol{\GG_1}$, where the first coefficient $a_1$ is given by
\beq\label{a_e1}
a_{1} = \frac{2}{(2\pi - \Ga){r_1 ^{\Gb}}} \int_{0}^{2\pi-\Ga}  \Gvf (Y(r_1,\Gt_1)) \sin \Gb \Gt_1 d\Gt_1
\eeq
for any $r_1 \in (0,1/2)$. We emphasize that $a_1$ is bounded independently of $\Ge$, which follows from \eqref{defGvf}  and Lemma \ref{lem_mini_v_E}. Since $\Bcal_1 (Y) = r_1^{\Gb} (Y) \sin \Gb \Gt_1 (Y) $, \eqnref{6000} can be rewritten as
\beq\label{6001}
(\Gvf -a_1 \Bcal_1) (Y) = \sum_{n=2}^{\infty} a_{n} r_1^{n\Gb} (Y) \sin n\Gb \Gt_1 (Y).
\eeq


For a given $\Gk < 1/ 2$, we choose $\Gk_1$ and $\Gk_2$ so that $\Gk <  \Gk_1  < \Gk_2  < 1/2$.
Lemma \ref{lem_mini_v_E} shows that there is $C_1>0$ such that
\beq\label{boundedness_of_U_inproof}
\left| \Gvf(Y) \right| \leq  |\p_2 \Ncal_{\Be}(Y)| +  |V (Y)| \leq C_1
\eeq
for all $Y \in B_{\Gk_2} (S_1) \setminus \ol{\GG_1} $.  Since $\Gvf = 0$ on $\p\GG_1$, there is a constant $\Gn$ such that for any $Y \in  \p  {B_{\Gk_1}(S_1) }  \setminus \ol{\GG_1}$, $\Gvf$ can be extended into $B_{\Gn} (Y)$ as a harmonic function. By shrinking $\Gn$ if necessary, we may assume that $\Gk_1+\Gn < \Gk_2$. Then the extended function is less than $2C_1$ in absolute value. Then, we obtain using \eqnref{standard} that
$$
|\nabla \Gvf (Y)| \leq C_2 := \frac{4C_1}{\eta} , \quad Y \in  \p B_{\Gk_1} (S_1)  \setminus \ol{ \GG_1}.
$$
We also have
$$
\left| a_{ 1}\nabla \Bcal_1(Y) \right| = \Gb \Gk_1 ^{\Gb-1}  |a_{ 1}|
   \leq   C_3 , \quad Y \in  \p B_{\Gk_1} (S_1)  \setminus \ol{ \GG_1}
$$
for some constant $C_3$. Thus, we have
\beq\label{the_boundedness_grd_U_B}
\norm{\nabla (\Gvf - a_1 \Bcal_1) }_{L^{\infty} \left( \p B_{\Gk_1} (S_1)  \setminus \ol{ \GG_1}\right)} \leq  C_2 + C_3 .
\eeq

Having \eqnref{6001} and \eqnref{the_boundedness_grd_U_B} in hand, we may adapt the same argument as in the proof of Theorem \ref{main2} (estimates in the region $R_2$) to infer that
there is a constant $C_4$ such that
$$
\left | \nabla (\Gvf - a_1 \Bcal_1)(Y) \right| \leq C_4, \quad Y \in  B_{\Gk} (S_1) \setminus \ol{\GG_1},
$$
since $\Gk  < \Gk_1$. Note that $|\nabla \Bcal_2(X)| \leq \Gb (1/2)^{\Gb-1}$ in $B_{\Gk} (S_1)$ since $\Gk < 1/2=|S_1-S_2|/2$. Thus we have
$$
\left | \nabla \Gvf (Y)  - a_{ 1} \nabla \left( \Bcal_1 (Y) + \Bcal_2 (Y) \right)\right| \lesssim 1 , \quad Y \in  B_{\Gk} (S_1) \setminus \ol{\GG_1}.
$$
We then have
$$
\left | \nabla \Gvf (Y)  - a_{ 1} \nabla \left( \Bcal_1 (Y) + \Bcal_2 (Y) \right)\right| \lesssim 1 , \quad Y \in  B_{\Gk} (S_2) \setminus \ol{\GG_2},
$$
thanks to symmetry of $\Gvf$ and $\Bcal_j$ with respect to $y$-axis. Thus \eqref{case3_variphi_B_1_B_2} follows.

We now prove that $a_1<0$ under the assumption (A). According to Lemma \ref {lem_mini_v_E}, $V$ attains its minimal value $\p_2 \Ncal_{\Be} (S_1)$ at $S_1$. By \eqref{defGvf} and the Taylor expansion of $\p_2 \Ncal_{ {\Be}}$ about $S_1$, we have
\begin{align}
\Gvf(Y) &=   \p_2 \Ncal_{ {\Be}} (Y) - V (Y) \nonumber \\
&\leq   \p_2 \Ncal_{ {\Be}} (Y) - \p_2 \Ncal_{ {\Be}} (S_1) \nonumber \\
& = \p_1\p_2 \Ncal_\Be(S_1) \left(y_1  - (- 1/2) \right) + \p_2 ^2 \Ncal_{ {\Be}} (S_1) y_2 + \Rcal_*(Y) , \label{GvfR*}
\end{align}
where the remainder $\Rcal_*(Y)$ satisfies
$$
\Rcal_*(Y)   \leq C_6 \left| Y-S_1 \right|^2
$$
for all $Y= (y_1,y_2) \in B_{1/ 3 } (S_1) \setminus \GG_1$ and  $C_6$ is a constant independent of small $\Ge>0$.

Since $\sin \Gb \Gt_1 > 0$ for $ \Gt_1 \in (0,\Ga) $, it follows from \eqref{a_e1} and \eqnref{GvfR*} that
\begin{align}
a_{1} &\le \frac{2}{(2\pi - \Ga){r_1^{\Gb}}} \Big[ \p_1\p_2 \Ncal_\Be(S_1) \int_{0}^{2\pi-\Ga}  \left(y_1(r_1,\Gt_1)  - (- 1/2) \right) \sin \Gb_1 \Gt_1 d\Gt_1 \nonumber \\
&\qquad + \p_2^2 \Ncal_\Be(S_1) \int_{0}^{2\pi-\Ga} y_2(r_1,\Gt_1) \sin \Gb_1 \Gt_1 d\Gt_1 + \int_{0}^{2\pi-\Ga} \Rcal_*(Y) \sin \Gb_1 \Gt_1 d\Gt_1 \Big] \nonumber \\
&:= \frac{2}{(2\pi - \Ga)^{-1}  {r_1^{\Gb}}} \big[ \p_1\p_2 \Ncal_\Be(S_1) I_1 +  \p_2^2 \Ncal_\Be(S_1) I_2 + I_3 \big]. \label{last}
\end{align}
We emphasize that the above inequality holds for all $r_1 <1/2$. One can immediately see from the oddness of the integrand that $I_2=0$. We also have
$$
|I_3| \le C_7 r_1^2 .
$$
On the other hand, we see that
$$
I_1= - r_1 \int_{0}^{2\pi - \Ga} \cos \left( \Gt_1 + \frac {\Ga} 2 \right)  \sin \left(\Gb \Gt_1 \right) d\Gt_1 = r_1 \frac  {2 \Gb}  {1-\Gb^2} \cos \frac{\Ga}{2} >0,
$$
and
$$
\p_1\p_2 \Ncal_\Be(S_1)= - \frac {8p} {\pi (1+4p^2)^2}< 0 .
$$
So, we have $I_1 > C_8 r_1$ for some positive constant $C_8$. We conclude from \eqnref{last} with a small enough $r_1$ that $a_1 <0$. This completes the proof.
\qed

\medskip
\noindent{\sl  Proof of Lemma \ref{lem:ACC}}.
By Lemma \ref{v_V_boundedness}, there exist constants $C_1$ and $\Ge_0$ such that
$$
\left| \Rcal_2 (X)\right|= \left| v(X) -  \frac 1 {\Ge}  V \left( {\Ge}^{-1}  X\right)\right| \leq  C_1
$$
for all $X  \in (\p B_{\mu}) \setminus (\GO_{1}\cup \GO_{2})$  and all $\Ge < \Ge_0$. We also have
$$
\left| \Rcal_2 (X)\right| = \left| v(X) -  \frac 1 {\Ge}  V \left( {\Ge}^{-1}  X\right)\right| = 0
$$
on $ (\p \GO_{1}\cup  \p \GO_{2}) \cap B_\mu $, where $\mu$ is the radius given in \eqref{GOone}. We then infer from the maximum principle that
\beq
\left| \Rcal_2 (X)\right| = \left| v(X) -  \frac 1 {\Ge}  V \left( {\Ge}^{-1}  X\right)\right| \leq  C_1 \label {rcal_2_2nd_bnded_L}
\eeq
for any $X  \in  B_\mu\setminus (\GO_{1}\cup \GO_{2})$  and any $\Ge < \Ge_0$.

We then use the same argument as for the proof of Lemma  \ref{estRtwo}, where an upper bound of $|\nabla \Rcal_2|$ is derived from an estimate of $|\Rcal_2|$. Here we take $\psi_{1+} = -\psi_{1-} = C_1$ on $\p B_{\mu} \setminus \overline {\GO_1}$ to infer that
$$
\left| \nabla \Rcal_2(X) \right| \lesssim \sum_{j=1}^2 \frac{1}{|X-V_j|^{1-\Gb}}, \quad X \in B_{1  } \setminus \ol {(\GO_{1}\cup \GO_{2})}.
$$
Thus \eqref{lem:ACC_ineq} follows.
\qed

 \end{document}